%% file: StochOpt-v2.tex
\documentclass[10pt,a4paper]{article}

\usepackage{hyperref}
\newcommand{\footremember}[2]{%
    \footnote{#2}
    \newcounter{#1}
    \setcounter{#1}{\value{footnote}}%
}
 
\input{StochOpt-v2_shared}

\input{zoom}

\date{}

\begin{document}

\title{Computational Aspects for Interface Identification Problems with Stochastic Modelling}
  
\author{Caroline Geiersbach\footremember{a}{University of Vienna, \texttt{caroline.geiersbach@univie.ac.at}} \and Estefania Loayza-Romero\footremember{b}{Chemnitz University of Technology, \texttt{estefania.loayza@math.tu-chemnitz.de}} \and Kathrin Welker\footremember{c}{Helmut-Schmidt-University / University of the Federal Armed Forces Hamburg, \texttt{welker@hsu-hh.de}}
}
\maketitle

\abstract{
In this paper, a shape optimization problem constrained by a random elliptic partial differential equation with a pure Neumann boundary is presented. The model is motivated by applications in interface identification, where we assume coefficients and inputs are subject to uncertainty. The problem is posed as a minimization of the expectation of a random objective functional depending on the uncertain parameters. A numerical method for iteratively solving the problem is presented, which is a generalization of the classical stochastic gradient method in shape spaces. Moreover, we perform numerical experiments, which demonstrate the effectiveness of the algorithm.
}

\bigskip


\bigskip


\section{Introduction}
Shape optimization is concerned with problems where an objective function is supposed to be minimized with respect to a shape, or a subset of $\R^d$. Many relevant problems contain additional constraints in the form of a partial differential equation (PDE), which describe the physical laws that the shape should obey. In applications, the material coefficients and external inputs might not be known exactly, but rather be randomly distributed according to a probability distribution obtained empirically. In this case, one might still wish to optimize over a set of these possibilities to obtain a more robust design. The random parameters are modeled with random fields, for example with a Karhunen-Lo\`{e}ve expansion \cite{Loeve1977, Loeve1978} or polynomial chaos \cite{Wiener1938}. The PDE constraint is then replaced by a set of PDE constraints, or a random PDE constraint, accounting for all possible realizations of the stochastic space.

One challenge in shape optimization is finding the correct model to describe the set of shapes; another is finding a way to handle the lack of vector structure of the shape space.
In principle, a finite dimensional optimization problem can be obtained for example by representing shapes as splines. However, this representation limits the admissible set of shapes, and the connection of shape calculus with infinite dimensional spaces \cite{Delfour-Zolesio-2001,Sokolowski1991} leads to a more flexible approach. Recently, it was suggested to embed shape optimization problems in the framework of optimization on shape spaces \cite{Schulz,Welker2016}. One possible approach is to cast the sets of shapes in a Riemannian viewpoint, where each shape is a point on an abstract manifold equipped with a notion of distances between shapes; see, e.g., \cite{MichorMumford2, MichorMumford1}. From a theoretical and computational point of view it is attractive to optimize in Riemannian shape manifolds because algorithmic ideas from \cite{Absil} can be combined with approaches from differential geometry. 
Here, the Riemannian shape gradient 
can be used to solve such shape optimization problems using the gradient descent method.
In the past, major effort in shape calculus has been devoted towards expressions for shape derivatives in the so-called Hadamard form, which are integrals over the surface (cf.~\cite{Delfour-Zolesio-2001,Sokolowski1991}). 
During the calculation of these expressions, volume shape derivative terms arise as an intermediate result. 
In general, additional regularity assumptions are necessary in order to transform the volume into surface forms.
Besides saving analytical effort, this makes volume expressions preferable to Hadamard forms.
In this paper, we consider the Steklov-Poincar\'{e} metric, which allows to use the volume formulations (cf.~\cite{SchulzSiebenbornWelker2015:2}).

Additional challenges arise in the stochastic setting. When the number of possible scenarios in the probability space is small, then the optimization problem can be solved over the entire set of scenarios. This approach is not relevant for most applications, as it becomes intractable if the random field has more than a few scenarios. For problems with random PDEs, either the stochastic space is discretized, or sampling methods are used. If the stochastic space is discretized, one typically relies on a finite-dimension assumption, where a truncated expansion is used as an approximation of the infinite-dimensional random field. Numerical methods include stochastic Galerkin method \cite{Babuska2004} and sparse-tensor discretization \cite{Schwab2011}. Sample-based approaches involve taking random or carefully chosen realizations of the input parameters; this includes Monte Carlo or quasi Monte Carlo methods and stochastic collocation \cite{Babuska2007}. In the stochastic approximation approach, dating back to a paper by Robbins and Monro \cite{Robbins1951}, one uses a stochastic gradient in place of a gradient to iteratively minimize the expected value over a random function. Recently, stochastic approximation methods have been proposed to efficiently solve PDE-constrained optimization problems involving uncertainty \cite{Geiersbach2019, Haber2012}. In this paper, we demonstrate a novel use of the stochastic gradient method on a model shape optimization problem. 

The paper is structured as follows. In~\Cref{sec:model_formulation}, a model problem will be presented; the model is an interface identification problem with uncertainty arising from PDE constraints containing random parameters and inputs. The problem will be formulated as a minimization problem over the expected deviation from a target measurement. The proposed algorithm for the numerical solution to the problem will be reviewed in~\Cref{sec:algorithmic_details}. Numerical experiments demonstrating the effectiveness of the method are shown in~\Cref{sec:numerical_experiments}. Finally, closing remarks are shared in~\Cref{sec:conclusion}.
\section{Model formulation}
\label{sec:model_formulation}

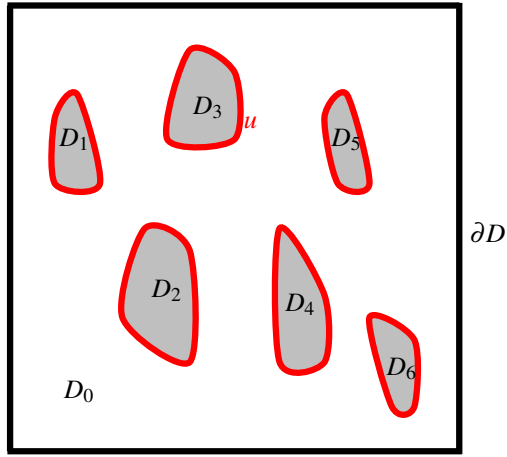
\begin{figure}
\begin{center}
   \begin{overpic}[width=.5\textwidth]{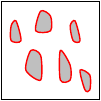}
   \put(12,68){$D_1$}
   \put(32,35){$D_2$}
   \put(41,75){$D_3$}
   \put(61,32){$D_4$}
   \put(71,68){$D_5$}
   \put(83,18){$D_6$}
   \put(13,13){$D_0$}
   \put(52,72){\textcolor{red}{$\Gi $}}
   \put(101,47){$\partial D $}
   \end{overpic}
   \end{center}
\caption{Illustration of the domain $D$}
\label{fig_domain}
\end{figure}

We consider a model interface identification problem, which has been studied in the deterministic setting in a number of texts \cite{Buttazzo2003,Ito-Kunisch-Peichl,Sokolowski1991}, and which we modify to allow for stochastic data.
We allow for a domain containing multiple shapes, leading to a multi-interface identification problem as illustrated in~\Cref{fig_domain}.
This model is for instance used in electrical impedance tomography, where the distribution of materials with different properties is examined based on measurements on the boundary.

Let $D\subset \mathbb{R}^2$ be a bounded domain with boundary $\partial D$. 
This domain is assumed to be Lipschitz and partitioned in an open subdomain $D_0\subset D$
and a finite number $N\in\mathbb{N}$ of disjoint open subdomains $D_i\subset D$ with variable boundaries $u_i:=\partial D_i$, $\sqcup_{i=0,\dots,N} D_i \bigsqcup\left(\sqcup_{i=1,\dots,N}u_i\right)=D$, where $\sqcup$ denotes the disjoint union.
For $i=1,\dots,N$, the union of all domains $D_i$ is denoted by $D_\text{int} := \sqcup_{i=1,\dots,N} D_i$ and the union of all boundaries $u_i$ is called the interface and denoted by $\Gi:=\sqcup_{i=1,\dots,N}u_i$.
\Cref{fig_domain} illustrates this situation.
In our setting, $D$ is meant to be composed of up to $N+1$ distinct materials. 
One should keep in mind that $D$ depends on $\Gi$, i.e., $D=D(\Gi)$. 
If $\Gi$ changes, then the subdomains $D_i\subset D$ change in a natural manner. 

For the model, we concentrate on one-dimensional smooth shapes.
In \cite{MichorMumford1}, the \emph{set of all one-dimensional smooth shapes} is characterized by
\begin{equation}
\label{B_e_2dim}
B_e= B_e(S^1,\R^2)\coloneqq \mathrm{Emb}(S^1,\R^2)/\mathrm{Diff}(S^1).
\end{equation}
Here, $\mathrm{Emb}(S^1,\R^2)$ denotes the set of all embeddings from the unit circle $S^1$ into $\R^2$, which contains all simple closed smooth curves in $\R^2$. Note that the boundary of the subdomain already characterize a shape. Thus, we can think of one-dimensional smooth shapes as images of simple closed smooth curves in the plane of the unit circle. In~\eqref{B_e_2dim}, $\mathrm{Diff}(S^1)$ is the set of all diffeomorphisms from $S^1$ into itself, which characterize all smooth reparametrizations. These equivalence classes are considered because we are only interested in the shape itself and images are not changed by reparametrizations. More precisely, shapes that have been translated represent the same shape. In contrast, shapes with different scaling are not equivalent in this shape space. 
In the following, we assume $u_i\in B_e$ for all $i=1,\dots,N$.

In this model, we allow different random variables accounting for source terms and material coefficients subject to uncertainty. It is assumed that probability distributions are known, for example by priorly obtained empirical samples. We allow for uncertainty in source terms and material constants by definition of a probability space $(\Omega, \mathcal{F}, \pP)$, where $\mathcal{F} \subset 2^{\Omega}$ is the $\sigma$-algebra of events and $\pP\colon \Omega \rightarrow [0,1]$ is a probability measure. 
In this paper, the probability space is to be understood as a product space $(\Omega, \mathcal{F}, \pP) = (\Omega_f \times \Omega_g \times \Omega_\kappa, \mathcal{F}_f \times \mathcal{F}_g \times \mathcal{F}_\kappa, \pP_{f} \times \pP_{g} \times \pP_{\kappa})$ as in~\cite{Gunzburger2014}. We assume without loss of generality that this product space is complete. 
We define $f\colon D \times \Omega_f \rightarrow \R$ and $g\colon\Go \times \Omega_g \rightarrow \R$ as a volume input and a boundary input function, respectively. We assume that these functions are finite dimensional noise in that there exist random vectors $\xi_f(\omega) = (\xi_f^1(\omega), \dots, \xi_f^{m_1}(\omega))$ with $\xi_f^i:\Omega \rightarrow \Xi_f^i \subset \R$ and $\xi_g(\omega) = (\xi_g^1(\omega), \dots, \xi_g^{m_2}(\omega))$ with $\xi_g^i:\Omega \rightarrow \Xi_g^i \subset \R$ such that
$$f(x,\omega) = f(x,\xi_f(\omega)), \quad  g(x,\omega) = g(x,\xi_g(\omega)) \quad \text{on } D\times \Omega.$$

On each subdomain $D_i$, the material is assumed to be constant, but randomly distributed. The function
\begin{equation}\label{eq:randomfield}
\kappa\colon D \times \Omega_\kappa \rightarrow \R,\  (x,\omega)\mapsto  \sum_{i=0}^N \kappa_{i}(\omega)\mathbbm{1}_{D_i}(x).
\end{equation}
defines the random material coefficient over the domain $D$, where $\kappa_i$ are independent real-valued random variables and $\mathbbm{1}_{D_i}$ denotes the indicator function of the set $D_i$. To simplify notation, we set $\xi:=(\xi_\kappa, \xi_f, \xi_g)$, $\Xi := \Xi_\kappa^{\text{in}} \times \Xi_\kappa^{\text{out}} \times \Xi_f^1 \times \cdots \times \Xi_f^{m_1} \times \Xi_g^1 \times \cdots \times \Xi_g^{m_2}$ and now write $\kappa(\xi) = \kappa(\cdot,\xi)$, $f(\xi) = f(\cdot,\xi)$, and $g(\xi) = g(\cdot,\xi)$ for a given $\xi \in \Xi$. 

For $J: B_e^N \times \Xi \rightarrow \R$, we use the notation $J(\cdot,\xi)\colon B_e^N \rightarrow \R$ to denote a deterministic functional for a single realization $\xi \in \Xi$. Here, $B_e^N$ stands for the cross-product of the shape manifold $B_e$ with itself $N$ times, one for each $u_i$. 
The expectation is then defined as the integral $\EE[J(u,\xi)]:=\int_\Omega J(u,\xi(\omega)) \, \d \pP(\omega).$
For the tracking-type objective functional
\begin{equation}
\label{objective}
J(u,\xi):=\frac{1}{2} \int_{D} (y(x,\xi)  - \bar{y}(x))^2 \, \dx
\end{equation}
we consider the following problem:
\begin{align}
&\underset{\text{with }u_i\in B_e}{\min_{\Gi=\sqcup_{i=1,\dots,N}u_i}} \quad \left\lbrace j(u):=\EE \big[ J(u,\xi) \big]\right\rbrace\label{eq:problem}
\\
&\hspace{.7cm}\text{s.t.} \quad y\colon D\times \Xi\to \R,\, (x,\xi)\mapsto y(x,\xi) \text{ satisfies} \qquad \nonumber
\end{align}
\begin{align}
- \nabla \cdot (\kappa \nabla y) &=  f,\quad \text{in } D \times \Xi \label{eq:PDE1} \\
\kappa \frac{\partial y}{\partial n} &= g, \quad \text{in } \partial D \times \Xi \label{eq:PDE2}
\end{align}
where $D=D(\Gi)$ depends on $\Gi$ and $\bar{y}\colon D\to\R$ denotes (deterministic) observed measurements. 
With the tracking-type objective functional $J(\cdot,\xi)$ the model is fitted to the data measurements $\bar{y}$.
The following continuity conditions are imposed explicitly for the state and flux at the interface:
\begin{equation}\label{eq:jumpconditions}
\left\llbracket \kappa \frac{\partial y}{\partial n} \right\rrbracket = 0, \quad \llbracket y \rrbracket  = 0, \quad  \text{in } \Gi \times \Xi,
\end{equation}
where the jump symbol $\left\llbracket\cdot\right\rrbracket$ is defined on $\Gi$ by $\llbracket v \rrbracket = v_1 - v_2$ with $v_1= \text{tr}_{1}(v)$ and $v_2 = \text{tr}_{2}(v)$ and $\text{tr}_1$, $\text{tr}_2$ are the trace operators defined on $u$ and $D_0$.

\begin{remark}
A regularization term such as perimeter regularization can be added to the objective functional~\eqref{objective}.
This is often necessary for analytical investigations, e.g., in order to establish that a problem is well-defined or to guarantee the existence of unique solutions.
In simulations, we observed that the scaling of the perimeter penalization needed to be set so small as to be negligible, so for simplicity we leave off this term.
\end{remark}

In order to solve model problem \eqref{eq:problem}--\eqref{eq:jumpconditions}, we need to calculate its shape derivative. 
We use the standard notation for the Sobolev space $H^1(D)$, where $H_0^1(D)$ indicates the subspace of $H^1(D)$ containing functions with disappearing trace, and its vector valued versions $H^1(D,\R^d),\,H_0^1(D,\R^d)$. For $H^r(D)$-valued random variables $y\colon H^r(D)\times\Xi \rightarrow \R$, we recall $y \in L^2(\Xi,H^r(D))$ if
$$\lVert y \rVert_{L^2(\Xi,H^r(D))}^2:= \mathbb{E}[\lVert y \rVert_{H^r(D)}^{r}] < \infty. $$
In addition, we define the Hilbert space  
$\hH := \{  v \in H^1(D) | \int_D v \, \d x= 0\}$.
In the following, $\xi \in \Xi$ is fixed but arbitrary. The weak formulation of the boundary value problem \eqref{eq:PDE1}--\eqref{eq:jumpconditions} for a fixed realization $\xi \in \Xi$ is: find $y = y(\cdot,\xi) \in \hH$ such that
\begin{equation}
\label{wf}
a_\xi(y,p)=b_\xi(p)\, , \ \forall p\in \hH
\end{equation}
with
\begin{align}
a_\xi(y,p) & = \int_D  \kappa(x,\xi) \nabla y(x,\xi)^T \nabla p(x) \, \dx   ,\label{bilinearform}\\
b_\xi(p) & = \int_D  f(x,\xi) p(x) \, \dx +\int_{\partial D} g(x,\xi) p(x) \, \ds ,\label{linearform}
\end{align}
where $\kappa$ is defined as in~\eqref{eq:randomfield}. 
In order for the shape derivative to be well-defined, we make the following technical assumptions: $\bar{y} \in H^1(D)$, $f \in L^2(\Xi,H^1(D))$ and $g \in L^2(\Xi,L^2(\partial D))$.

The Lagrangian of~\eqref{eq:problem}--\eqref{eq:jumpconditions} is defined as
\begin{equation}
 \LL_\xi(u,y,p):= J(u,\xi)+a_\xi(y,p)-b_\xi(p),
\label{lagrangian_obj}
\end{equation}
where $J(\cdot,\xi)$ is defined in~\eqref{objective}, $a_\xi$ in~\eqref{bilinearform} and $b_\xi$ in~\eqref{linearform}.
For any $\xi \in \Xi$, a saddle point $(y,p)\in \hH\times \hH$ of the Lagrangian is given by
\begin{align}
\frac{\partial \LL_\xi(u,y,p)}{\partial y}=\frac{\partial \LL_\xi(u,y,p)}{\partial p}=0,\label{saddlepointcond}
\end{align}
which leads to the \emph{state equation}~\eqref{eq:PDE1}--\eqref{eq:jumpconditions} and the \emph{adjoint equation}
\begin{align}
-   \nabla \cdot (\kappa \nabla p) &=  \bar{y} -y,\quad \text{in }D \times \Xi \label{eq:adjointPDE1} \\
\kappa \frac{\partial p}{\partial n}  &= 0,\hspace{.5cm} \quad\text{in } \Go \times \Xi \label{eq:adjointPDE2}
\end{align}
for $p\colon D\times \Xi\to \R,\, (x,\xi)\mapsto p(x,\xi)$.
As with the state equation, we have interface conditions for the adjoint equation:
\begin{equation}\label{eq:jumpconditions-adjoint}
\left \llbracket \kappa \frac{\partial p}{\partial n} \right\rrbracket = 0, \quad \llbracket p \rrbracket  = 0,\quad \text{in } \Gi \times \Xi.
\end{equation}
The \emph{design equation} is given by the shape derivative. If we consider the perturbation of identity, the shape derivative of $\LL$ at $u$ and fixed $\xi \in \Xi$ in the direction of a sufficiently smooth vector field $V$ is defined by
	\begin{equation}\label{eq:Eulerianderivative}
	d\LL_\xi(u,y,p)[V]:= \lim_{t \rightarrow 0^+} \frac{\LL_\xi(F_t(u),y,p) - \LL_\xi(u,y,p)}{t},
	\end{equation}
where $F_t(u) := \{\text{id}(x)+tV(x)\colon x\in u\}$ and $F_0(u)=u$.
Using standard techniques for calculating the shape derivative (min-max approach~\cite{Delfour-Zolesio-2001}, chain rule approach~\cite{SokoZol}, rearrangement method~\cite{Ito-Kunisch-Peichl}), we get the shape derivative of our model problem in volume formulation:
\begin{equation}
\label{sd_j1}
\begin{split}
d J(u,\xi)[V]=& \int_D -\kappa\nabla y^T(\nabla V+\nabla V^T)\nabla p - (y-\bar{y})\nabla\bar{y}^TV-\nabla f^TV p \\
& \hspace{.4cm}+ \text{div}(V)\left( \frac{1}{2} (y-\bar{y})^2 + \kappa\nabla y^T\nabla p - fp \right)\,\dx.
\end{split}
\end{equation}

\begin{remark}
Note that the shape derivative arises in two equivalent notational forms:
\vspace{.1cm}
\begin{align*}
d J(u,\xi)[V]&\coloneqq \int_{D} R(x)V(x)\, \dx  &\text{(volume/weak formulation),}\\[5pt]
d J(u,\xi)[V]&\coloneqq\int_{\Gi} r(s)\left<V(s), n(s)\right> \ds  &\text{(surface/strong formulation).}
\end{align*}
Here, $R$ is a differential operator acting linearly on the vector field $V$ and  $r\in L^1(\Gi)$.
In this paper, we use the volume formulation.
\end{remark}

\begin{remark}
\label{remark:perturbation_of_identity}
In this paper, we consider the perturbation of identity. Of course, it is also possible to use, e.g., the speed method in order to transform the domain $D$. In this case, the definition of $F_t$ changes.
\end{remark}

\begin{remark}
As mentioned in the introduction, the algorithm proposed in this work is basically a generalization of the stochastic gradient method. For this method, descent directions can be chosen as the gradient of $J(\cdot,\xi)$ for a randomly chosen $\xi$. This will be clarified in the next section.
This is the reason why the shape derivative is computed just for $J(\cdot,\xi)$ and not for $j$.
\end{remark}

Now that we have the state, adjoint and design equations, we are ready to solve our problem computationally. In the next section, we will present our proposed method for the computational solution of the problem, including the discretization scheme, optimization method, and post-processing techniques for the improvement of mesh aspect ratios.
\section{Algorithmic and computational details}
\label{sec:algorithmic_details}
This section is split into two subsections. \Cref{subsection_Algorithm} is devoted to the presentation of the algorithmic aspects for the solution of the problem described in~\cref{sec:model_formulation}. To our knowledge, this is a novel application of the stochastic gradient method to a shape optimization problem under uncertainty. 
In~\cref{subsection_CompAspects}, we focus on computational aspects and combine these with the algorithmic details of~\cref{subsection_Algorithm} resulting in a stochastic shape gradient method for a discretized version of the model problem.

\subsection{Algorithmic aspects}
\label{subsection_Algorithm}

The \emph{stochastic gradient method} is an algorithm for solving problems of the form 
\begin{equation}
\label{eq:stochastic-optimization-problem}
\min_{x}\,\EE[J(x,\xi)].
\end{equation}
The method uses iterations of the form $x^{n+1} = x^n - t^n \nabla J(x^n,\xi^n)$, where at iteration~$n$ the \textit{stochastic gradient} at $x^n$ for a random sample ${\xi^n} \in \Omega$, is given by $\nabla J(x^n,\xi^n) \approx \nabla \EE[J(x^n,\xi)]$. We make the assumption that it is possible to generate independent and identically distributed (i.i.d.) samples $\xi^n$ from a known probability distribution.
This method dates back to~\cite{Robbins1951} and has been used in many different applications, notably as well in optimization on manifolds~\cite{Bonnabel2011} and in PDE-constrained optimization~\cite{Geiersbach2019}. Much research has been devoted to the proper choice of the step size $t^n$. The least involved choice is the so-called Robbins-Monro step size, 
\begin{equation}
\label{eq:Robbins-Monro}
t^n \geq 0, \quad \sum_{n=0}^\infty t^n = \infty, \quad \sum_{n=0}^\infty ({t^n})^2 < \infty,
\end{equation}
which was also used in the original work~\cite{Robbins1951} to establish almost sure convergence of the sequence $\{x^n \}$ to a stationary point for~\eqref{eq:stochastic-optimization-problem}.  It is notable that the Armijo rule $t^n:=\alpha \rho^{m_n}$ with minimal $m_n \in \mathbb{N}_0$ satisfying
\begin{equation}
\label{eq:armijo-rule}
J(x^n-t^n \nabla J(x^n,\xi^n), \xi^n) \leq J(x^n,\xi^n) - t^n c \lVert \nabla J(x^n,\xi^n) \rVert^2
\end{equation}
with $\alpha > 0$, $\rho \in (0,1)$, $c \in (0,1)$ may not converge without additional variance reduction techniques. If instead of choosing a single realization at iteration $n$, one chooses $N_n$ samples $\xi^n = (\xi^{n,1}, \dots, \xi^{n,N_n})$ with $N^n \rightarrow \infty$ as $n \rightarrow \infty$, then convergence using~\eqref{eq:armijo-rule} can be demonstrated under relatively weak assumptions \cite{Wardi1990,Shapiro1996}. Note that this involves increasingly expensive calculations for the estimate of the objective function $J(x,\xi^n)= 1/N_n\sum_{l=1}^{N_n} J(x, \xi^{n,l})$ and likewise for the stochastic gradient $\nabla J(x,\xi^n).$ In numerical simulations, we will compare both of these approaches. Further step size rules can be reviewed in a nice summary by George and Powell~\cite{George2006}.
 
Now, we formalize what is meant by stochastic gradient for our problem. In this work, we focus on one-dimensional smooth shapes which are elements of the shape space $B_e$ defined in~\eqref{B_e_2dim}. This shape space is a Riemannian manifold (cf.~\cite{MichorMumford1,MichorMumford}).
If we want to optimize on a Riemannian shape manifold, we have to find a representation of the shape derivative with respect to the Riemannian metric under consideration, called the \emph{Riemannian shape gradient},
which is required to formulate optimization methods on a shape manifold.
In~\cite{SchulzSiebenbornWelker2015:2}, the authors present a metric based on the Steklov-Poincar\'{e} operator, which allows for the computation of the Riemannian shape gradient as a representative of the shape derivative in volume form.
Besides saving analytical effort during the calculation process of the shape derivative, this technique is computationally more efficient than using an approach which needs the surface shape derivative form. For example, the volume form allows us to optimize directly over the hold-all domain $D$ containing one or more elements $u_i \in B_e$, whereas the surface formulation would give us descent directions (in normal directions) for the boundary $u_i$ only, which would not help us to move mesh elements around the shape. Additionally,  when we are working with a surface shape derivative, we need to solve another PDE in order to get a mesh deformation in the hold-all domain $D$.
Following the ideas presented in~\cite{SchulzSiebenbornWelker2015:2}, we choose the Steklov-Poincar\'e metric, denoted by $G^S$.
In the setting of the shape space $B_e$, the mesh deformation vector $V\in H_0^1(D,\R^2)$ can be viewed as an extension of a Riemannian shape gradient to the hold-all domain $D$ because of the identities 
\begin{equation}
\label{Steklov_identity}
G^S(v,u)= d J(u,\xi)[U]=a(V,U)\quad \forall U\in H_0^1(D,\R^2),
\end{equation}
where $v=(\text{tr}(V))^Tn,u=(\text{tr}(U))^Tn$ with $\text{tr}(\cdot)$ denotes again the trace operator on the Sobolev spaces.
Here, $a(\cdot,\cdot):H_0^1(D,\R^2) \times H_0^1(D,\R^2) \rightarrow \R$ is a symmetric and coercive bilinear form.
One option for the operator $a(\cdot, \cdot)$ is chosen to be the bilinear form associated with the linear elasticity problem, i.e.,
\begin{equation}
\label{eq:shape_gradient}
a^{\text{elas}}(V, U):=\int_D (\lambda \text{tr}  ( \epsilon(V)  )\id + 2  \mu \epsilon(V)) : \epsilon(U)\, \dx,
\end{equation}
where $\epsilon(U)\coloneqq \frac{1}{2} \, (\nabla U + \nabla U^T)$, $A : B$ denotes the Frobenius inner product for two matrices $A, B$ and $\lambda,\mu$ denote the Lam\'{e} parameters.
To summarize, for a fixed $\xi \in \Xi$, we need to solve the following so-called \emph{deformation equation}: find $V \in H_0^1(D,\R^2)$ such that
\begin{equation}
a^{\text{elas}}(V, U) = d J(u,\xi)[U] \quad \forall U\in H_0^1(D,\R^2).
\label{deformatio_equation}
\end{equation}

The main advantage of this Steklov-Poincar\'{e} metric approach is that the identity~\eqref{Steklov_identity} holds, meaning that the Riemannian metric $G^S(\cdot, \cdot)$, which is naturally defined over the interfaces, can be equivalently reformulated in terms of the bilinear form $a(\cdot, \cdot)$ over the whole domain. This last observation is the main approach we will use in the numerical solution of our model problem.

\begin{remark}
\label{rmk:exponential_retractions}
In general, we need the concept of the exponential map and vector transports in order to formulate optimization methods on a shape manifold. The calculations of optimization methods have to be performed in tangent spaces
because manifolds are not necessarily linear spaces. This means points from
a tangent space have to be mapped to the manifold in order to get a new
shape-iterate, which can be performed with the help of the exponential map. However, the computation of the exponential map is prohibitively expensive in the most applications because a calculus of variations problem must be solved or the Christoffel symbols need to be known.
It is much easier and much faster to use a first-order approximation of the exponential map.
In \cite{Absil}, it is shown that a so-called \emph{retraction} is such a first-order approximation and sufficient in most applications. 
We refer to \cite{SchulzWelker}, where a suitable retraction on $B_e$ is given.
\end{remark}

\subsection{Computational aspects}
\label{subsection_CompAspects}

In order to numerically solve problem~\eqref{eq:problem}--\eqref{eq:jumpconditions}, we need to consider some kind of discretization; we use the finite element method (FEM). Let $D_h$ be a triangulation of the domain $D$ with a partition of subsets $D_{i,h}:=\overline{D_i} \cap D_h$ such that $D_h = \sqcup_{i=0,\dots,N} D_{i,h}$. Moreover, let $u_{i,h}:= {D_{i,h}} \cap {D_{0,h}}$ denote the discretization of the interfaces for $i=1, \dots, N$ and $\partial D_h := \partial D \cap D_h$. We assign $D_{\text{int},h}:=\sqcup_{i=1, \dots, N} D_{i,h}$ and $u_{h}: = \sqcup_{i=1, \dots, N} u_{i, h}.$
We consider the set, 
\[
W_h\coloneqq \left\{ y \in H^1(D) \, \big| \,  y|_T \in \mathcal{P}_1 \text{ for all  } T \text{ in } D_h\right\}.
\]
where $\mathcal{P}_1$ is the space of all first degree polynomials. Furthermore, let us define the spaces 
\[
\hHh \coloneqq \left\{ y \in \hH \, \big| \,  y|_T \in \mathcal{P}_1 \text{ for all  } T \text{ in } D_h\right\}
\] 
where $\hH := \{  v \in H^1(D) | \int_D v \, \d x= 0\}$, and
\[
W_{0,h} \coloneqq \left\{ V \in H_0^1({D, \R^2}) \, \big| \,  V_i|_T \in \mathcal{P}_1 \text{ for all } i=1,2 \text{ and } T \text{ in } D_h\right\}.
\]
Once we have discretized the domain and found an approximation for the Sobolev spaces, we are able to present the discrete version of the state equation, which is given in its weak formulation as follows: 
\begin{equation}
\label{weak_state}
\text{Find } y_h\in \hHh\colon \int_{D_h} \kappa_h \nabla y_h^T \nabla p_h \, \dx = \int_{D_h} f_h p_h\, \dx + \int_{\partial D_h} g_h p_h \, \ds \quad\forall p_h \in \hHh,
\end{equation}
where $\mathbb{I}_h: C(\bar{D}) \rightarrow W_h$ is the piecewise linear interpolation of a continuous function in the space $W_h$ and $\kappa_h := \mathbb{I}_h(\kappa)$, $f_h := \mathbb{I}_h(f)$, $g_h := \mathbb{I}_h(g)$.
In the same way, the adjoint state will be computed by solving the problem
\begin{equation}
\label{weak_adjoint}
\text{Find } p_h\in \hHh \colon \int_{D_h} \kappa_h \nabla p_h^T \nabla v_h \, \dx = \int_{D_h} (y_h - \bar{y}_h) v_h\,\dx, \quad\forall v_h \in \hHh
\end{equation}
where $\bar{y}_h:=\mathbb{I}_h(\bar{y})$.
Furthermore, as explained in the previous section, these equations need to be solved for a random $\xi$. For readability, we have suppressed the dependence of $\kappa_h$,  $y_h$, and $p_h$ on $\xi$.

Once we have computed the state and adjoint variables, the shape derivative will be given by the discrete counterpart of equation~\eqref{sd_j1}, i.e., 
\begin{equation}
\label{discrete_sd_j1}
\begin{split}
dJ(u_h,\xi)[U_h]=& \int\limits_{D_h} -\kappa_h\nabla y_h^T(\nabla U_h+\nabla U_h^T)\nabla p_h - (y_h-\bar{y}_h)\nabla\bar{y}_h^TU_h-\nabla f_h^TU_h p_h \\
& \hspace{.4cm}+ \text{div}(U_h)\left( \frac{1}{2} (y_h-\bar{y}_h)^2 + \kappa_h\nabla y_h^T\nabla p_h - f_h p_h \right)\,\dx
\end{split}
\end{equation}

Finally, the discretized version of the deformation equation is given by: find $V_h \in W_{0,h}$ such that
\begin{equation}
\label{discretized_deformation_equation}
a^\text{elas}(V_h, U_h)= dJ(u_h,\xi)[U_h]\quad\forall U_h\in W_{0,h}.
\end{equation}
As mentioned in Remark~\ref{rmk:exponential_retractions}, after the computation of the deformation field, we need to use the exponential function or an approximation of it (like a retraction mapping) in order to update the shape for the next iteration. In this context, following the ideas of~\cite{SchulzWelker}, we have chosen a retraction that in its discretized version acts like the perturbation of identity method.
That means we assign a direction given by $V_h$ to each of the nodes of the triangulation $D_h$ and move them with certain step size. 
In order to guarantee good behavior of the method (invertibility of the operator $F_t$ from the perturbation by identity), this step sizes need to be small; otherwise, the algorithm could produce low aspect ratios or overlapping elements. One of the principal reasons for this behavior is that the perturbation by identity method does not consider mesh connectivity, or how close a node is to its neighbors.
In order to overcome this drawback, we use the following techniques.  

As discussed in~\cite{SchulzSiebenbornWelker2015:2}, an unmodified right-hand side of the discretized deformation equation leads to deformation fields causing meshes with bad aspect ratios. Thus, we set values of the discretized shape derivative to zero if the corresponding element does not intersect with the interface, i.e.,
$$dJ(u_h,\xi)[U_h]=0\quad \forall U_h\text{ with }\text{supp}(U_h)\cap u_{h}=\emptyset.$$
Furthermore, in order to ensure good behavior of the mesh along iterations of the algorithm, we  follow ideas from \cite{SchulzSiebenborn2016}; for each iteration step $n$, we choose the Lam\'{e} parameters as follows: $\lambda=0$
and $\mu \in [\mu_{\min},\mu_{\max}]$ decreasing smoothly from $\Gi$
to the outer boundary. One possible way to model this behavior is to solve the following Poisson equation
  \begin{align}
  \nonumber\Delta \mu_h &= 0\hspace{.9cm}   \text{ in } D_{0,h} \sqcup D_{\text{int},h}\\
  \label{eq:system_mu_elas}\mu_h &= \mu_{\max}  \hspace{0.4cm}\text{ on } u_{h}\\
  \nonumber\mu_h &= \mu_{\min} \hspace{0.5cm} \text{ on } \partial D_h.
  \end{align} 
In~\cref{subsection_InfluenceLame}, we will present an experiment that shows the importance of the choice for $\mu_{\min}$ and $\mu_{\max}$. Correct choices for these parameters ensure that the cells in the mesh maintain good aspect ratios throughout the optimization process.

\begin{algorithm}
	\begin{algorithmic}[0]
		\State \textbf{Initialization:} Choose $D^0=\sqcup_{i=0,\ldots,N}D_{i,h}^0\sqcup_{i=1,\ldots,N} u_{i,h }^0$ 
		\For{$n=0,1,\ldots$}
		\State Generate $\xi^n\in  \Xi$, choose $t^n \geq 0$
		\State $y_h^n \gets$ Solve the state equation given in~\eqref{weak_state} with $\xi = \xi^n$
		\State $p_h^n \gets$ Solve adjoint random equation given in~\eqref{weak_adjoint} using $y_h = y_h^n$ and $\xi = \xi^n$
		\State $dJ_{\xi_n}(D_h^n)[U_h] \gets$ Assemble the random discrete shape derivative \eqref{discrete_sd_j1} with $U_h$ such that \\
		\hspace{25mm} $\text{supp}(U_h)\cap u_{h}=\emptyset$ 
		\State $\mu_h^n \gets$ Solve Poisson equation given in~\eqref{eq:system_mu_elas}
		\State $V_h^n \gets$ Compute the mesh deformation vector field by solving the deformation equation~\eqref{discretized_deformation_equation}
		\State Set $D_h^{n+1} = D_h^n - t^n V_h^n$
\EndFor
\end{algorithmic}
\caption{Stochastic shape gradient method for the model problem}
\label{alg:stochastic_gradient_shape_optimization}
\end{algorithm}

In~\cref{alg:stochastic_gradient_shape_optimization}, we present the \emph{Stochastic Shape Gradient Method (SSGM)} formulated specifically for the discrete version of the model problem described in~\cref{sec:model_formulation}.

The model problem (\ref{eq:problem})-(\ref{eq:PDE2}) is formulated for one-dimensional shapes $u_i\in B_e$, $i=1,\dots,N$. This means that we are dealing with objects in 2D. Of course, it is possible to generalize our problem to 3D shapes. For a generalization, we refer to to~\cite{Siebenborn2017}, which considers a multiple shape interface problem in 3D. One focus of this paper is the performance of 3D computations. Here, the effectiveness of the computations is guaranteed by using multigrid methods and parallelization techniques. We want to mention that our proposed method can also be applied to higher-dimensional objects. However, one needs to take into account that problems occur, which are due to the high dimension. In 3D, we have to deal with many degrees of freedom. Thus, our method needs some tuning in higher dimensions to solve the problem efficiently. For example, parallelization techniques or multigrid methods can be performed, as suggested by~\cite{Siebenborn2017}. Since addressing this kind of question goes beyond the scope of this paper, we restrict ourselves to 2D.
Although the stochastic gradient, the numerical solution of PDEs, constrained shape optimization problems had been widely studied in the past, there have not been a combination of both, and the scientific computing is challenging as the experiments in the next section show. We want to highlight that the main advantage of our method is the lower computational cost added by considering random parameters. This fact is obtained due to the low complexity of the stochastic gradient. Furthermore,  up-to-the authors' knowledge, there has not been an extensive study on the different techniques for preventing mesh destruction and bad aspect ratios like the ones presented here.

\begin{remark}
We want to mention that model problem (\ref{eq:problem})-(\ref{eq:PDE2}) for $N=1$ in higher dimensions is investigated analytically in \cite{GeiersbachLoayzaWelker}. In particular, an asymptotic convergence result on Riemannian manifolds is proven.
\end{remark}

\section{Numerical experiments}
\label{sec:numerical_experiments}
The principal aim of this section is to 
show the behavior of~\cref{alg:stochastic_gradient_shape_optimization} through numerical experiments. This section is divided into five experiments. The first one shows the applicability to a problem with multiple shapes. The second one presents an analysis of the behavior of the meshes for different values of the parameters $\mu_{\min}$ and $\mu_{\max}$ in the equation~\eqref{eq:system_mu_elas}. The third shows convergence to different target shapes. The fourth shows the influence of the individual random variables on the outcome. Finally, we look at results obtained using different step size rules and probability distributions with high variance.

Throughout this section we choose $D=(0,1)^2$ and discretize as described in~\cref{subsection_CompAspects} with up to 10,000 triangles. In all experiments, the volume input $f$ is set to $f \equiv 0.$ We use $\mathcal{N}(\rho,\sigma,a,b)$ to denote a truncated normal distribution with mean $\rho$, standard deviation $\sigma$, and lower and upper bounds denoted by $a$, $b$, respectively. We consider a simplified version of the model where $\kappa_{\text{int}}:=\kappa_1= \cdots = \kappa_N.$ The measurement $\bar{y}$ is generated by solving the state equation~\eqref{weak_state} on a target shape $\bar{u}$ with constant parameters $\kappa_0 = 1.5$, $\kappa_{\text{int}} = 4$ and $g= 10$. Unless otherwise indicated, the bounds of the Lam\'{e} parameter $\mu$ were set to $\mu_{\min} = 10 $ and $\mu_{\max}= 25$. 

In all experiments, we use a single realization $\xi^n$ at each iteration $n$ to generate the stochastic shape gradient. For some convergence plots, we approximate function values $\EE[J(u^n, \xi)]$ and the $L^2$-norm of the shape gradient $\EE[\lVert V^n \rVert_{L^2(D,\R^2)}]$ using additional sampling at each point $u^n$ by the estimates
$$\hat{j}_n :=\frac{1}{m} \sum_{l=1}^m J(u^n,\xi^{n,l}) \approx \mathbb{E}[J(u^n,\xi)], \quad\hat{v}_n:=\frac{1}{m} \sum_{l=1}^m \lVert V_n(\xi^{n,l)}) \rVert_{L^2(D,\R^2)} \approx \mathbb{E}[\lVert V^n \rVert_{L^2(D,\R^2)}]$$
with $m$ i.i.d.~samples $\{\xi^{n,1}, \dots, \xi^{n,m}\}$ generated at each $n$.

\subsection{Multiple shapes}

\begin{figure}
\centering
\includegraphics[scale=0.25]{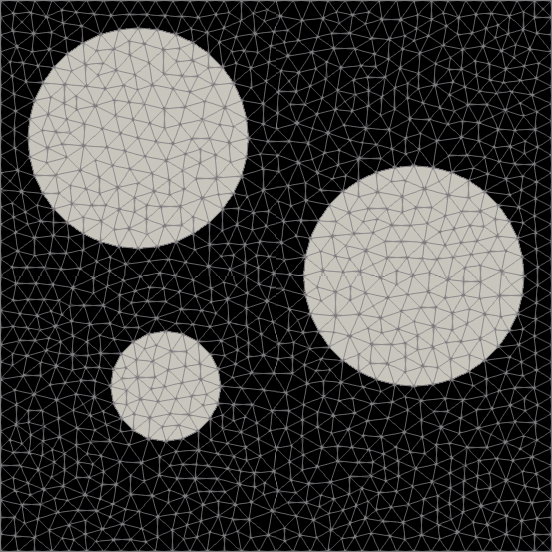}
\caption{Target shape $\bar{u}$}
\label{fig:target_experiment1}
\end{figure}

The main objective of this experiment is to present the numerical solution of problem~\eqref{eq:problem}-~\eqref{eq:jumpconditions} for the target configuration depicted in~\Cref{fig:target_experiment1}. 
Additionally, we present experiments on meshes with different sizes
to observe their effect on the solution. We used the Armijo line search rule~\eqref{eq:armijo-rule} with $\alpha=50$, $\rho = 0.5$, and $c= 10^{-4}$ and the random parameters $\kappa_0  \sim \mathcal{N}(1.5,10^{-2},1,2) $, $ \kappa_{\text{int}}\sim \mathcal{N}(4,10^{-2},3,5)$ and $g \sim \mathcal{N}(10,10^{-2},9,11)$.

\begin{figure}
\centering
\begin{subfigure}[b]{0.3\textwidth}
\includegraphics[scale=0.18]{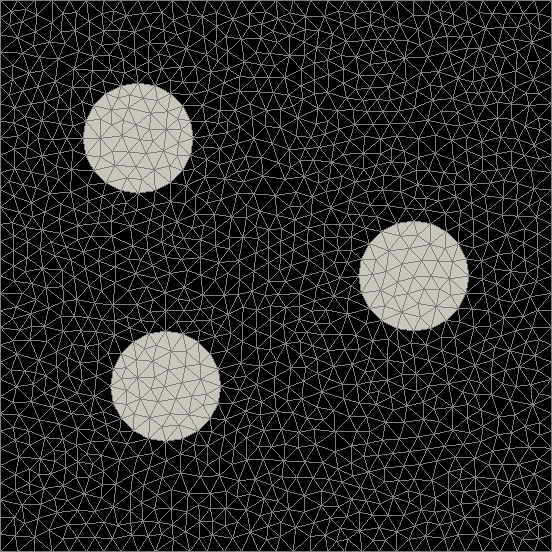}   
\caption{$D^0$} 
\label{fig:initial_experiment1}
\end{subfigure}
~
\begin{subfigure}[b]{0.3\textwidth}
\includegraphics[scale=0.18]{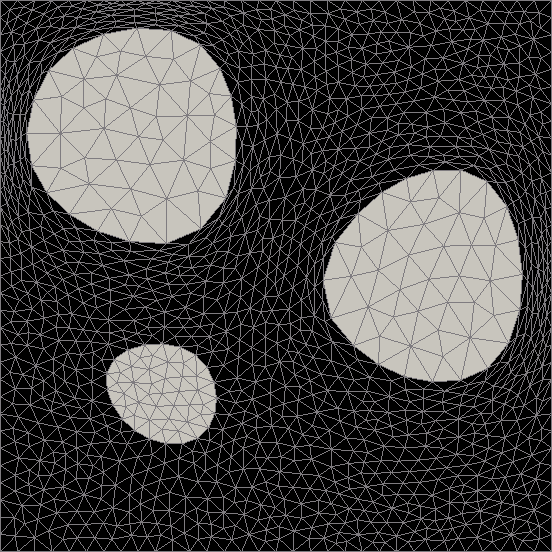}   
\caption{$D^{150}$} 
\end{subfigure}
~
\begin{subfigure}[b]{0.3\textwidth}
\includegraphics[scale=0.18]{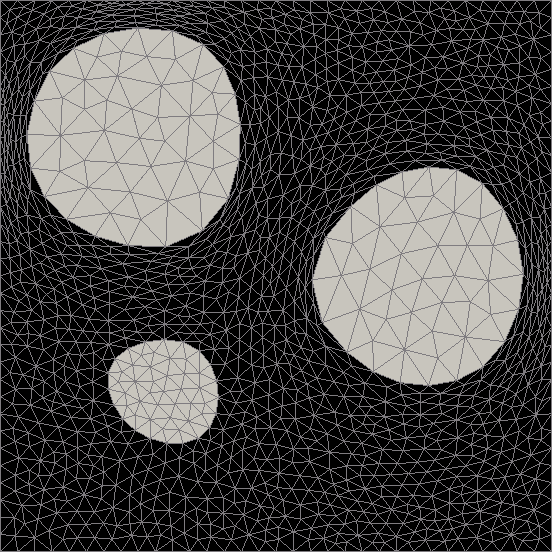}
\caption{$D^{300}$}
\label{fig:final_experiment1}
\end{subfigure}
\caption{Shapes obtained on a coarse mesh}
\label{fig:experiment1}
\end{figure}

\begin{figure}
\begin{subfigure}[b]{0.45\textwidth}
\includegraphics[scale=0.45]{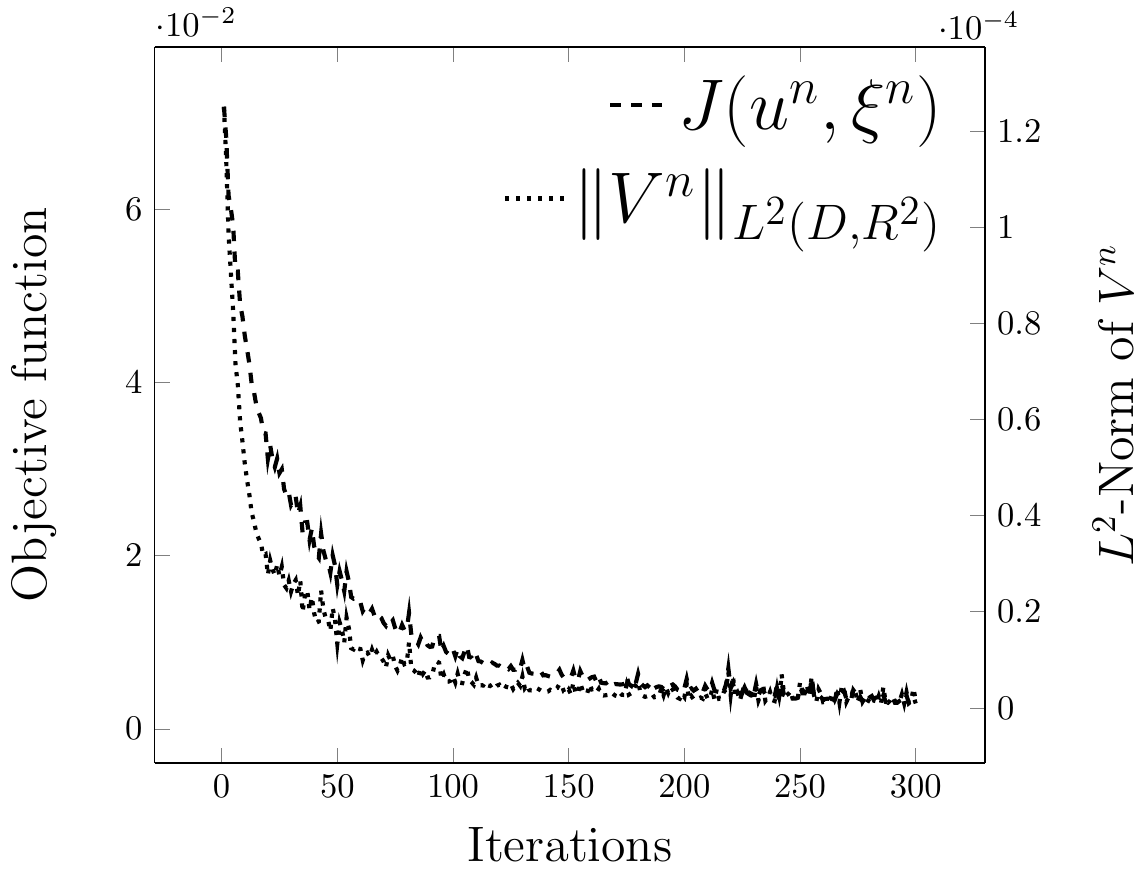}
\caption{Mesh with 3000 elements}
\label{fig:decay}
\end{subfigure}
~
\begin{subfigure}[b]{0.45\textwidth}
\includegraphics[scale=0.45]{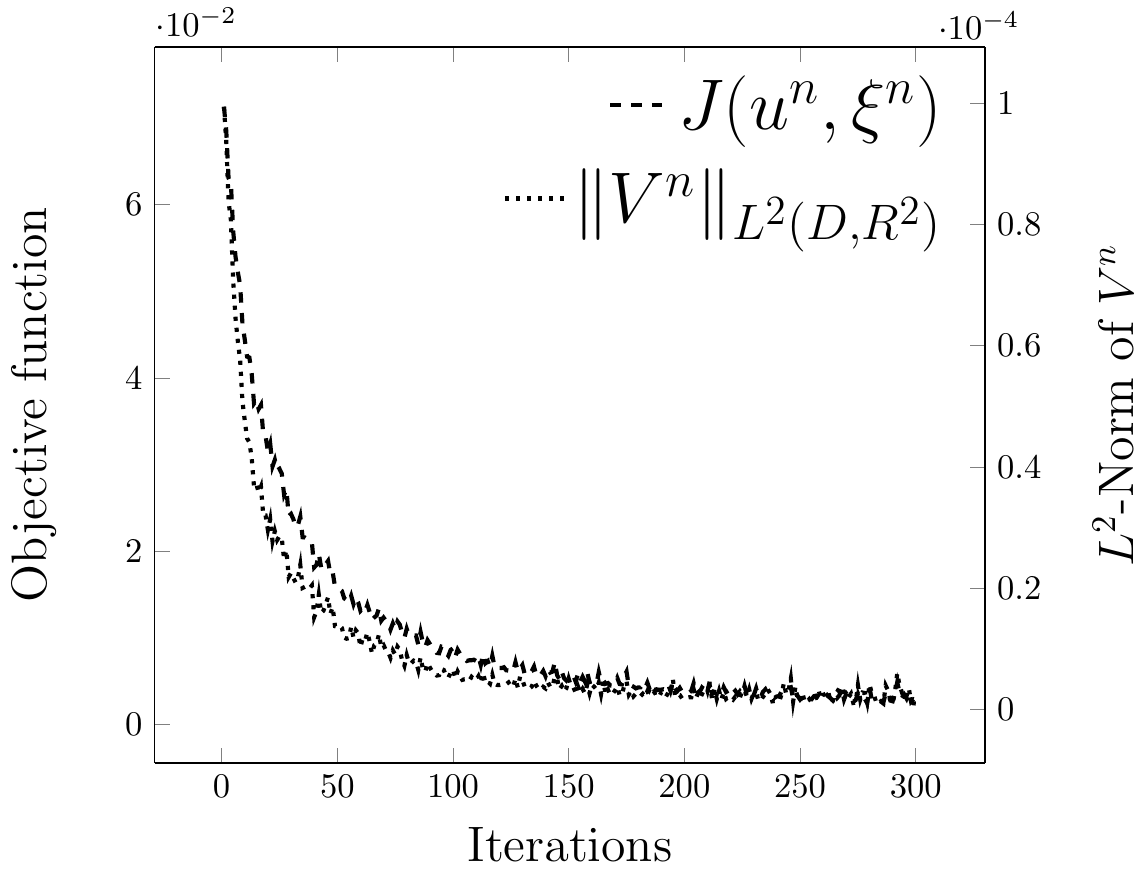}
\caption{Mesh with 10,000 elements}
\label{fig:decay_fine}
\end{subfigure}
\caption{Convergence plots}
\label{fig:decay_experiment1}
\end{figure}

First, we used a mesh with approximately 3,000 elements. We let the algorithm iterate 300 times, and~\Cref{fig:experiment1} shows the initial, an intermediate and the final shapes obtained. A convergence plot for a single sample of the objective function value and stochastic shape gradient is shown in~\Cref{fig:decay}. Furthermore, we use $m=100$ samples to approximate the objective function at the last iteration to find $\hat{j}_{300} \approx 0.00323.$

\begin{figure}
\centering
\begin{subfigure}[b]{0.3\textwidth}
\includegraphics[scale=0.18]{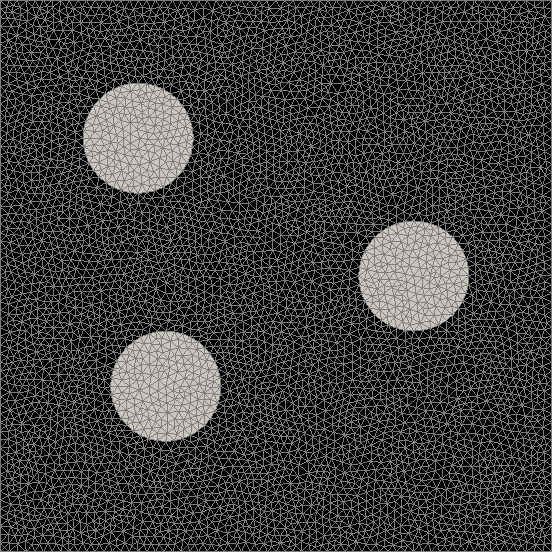}   
\caption{$D^0$} 
\end{subfigure}
~
\begin{subfigure}[b]{0.3\textwidth}
\includegraphics[scale=0.18]{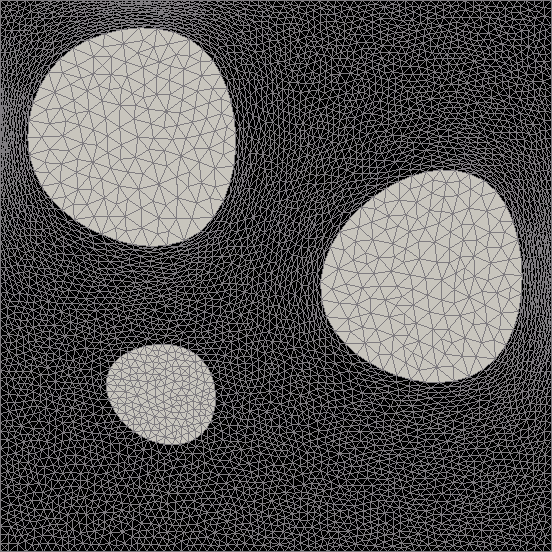}
\caption{$D^{150}$}
\end{subfigure}
~
\begin{subfigure}[b]{0.3\textwidth}
\includegraphics[scale=0.18]{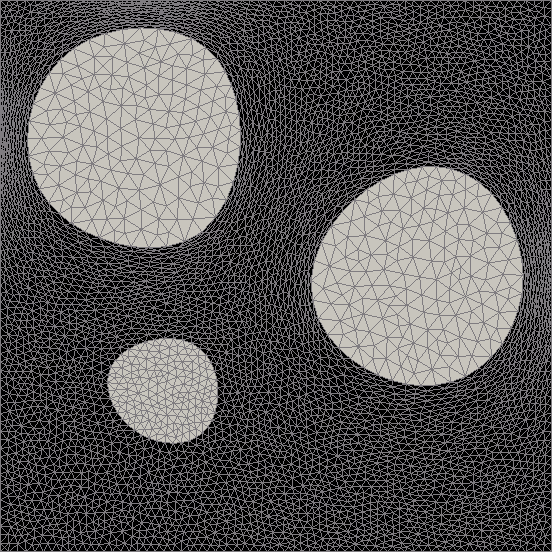}
\caption{$D^{300}$}
\end{subfigure}
\caption{Shapes obtained on a fine mesh}
\label{fig:experiment1_fine}
\end{figure}

Second, we would like to show that the solutions can be improved just by taking a finer mesh for both the initial and target shapes, in this case with around 10,000 elements. In~\Cref{fig:experiment1_fine}, we present the obtained solutions. Moreover, the convergence plot for this experiment is showed in~\Cref{fig:decay_fine}, and the value of the objective function is approximated with $m=100$ samples with $\hat{j}_{300} \approx 0.00311$. As expected, we observe a lower objective function value with the finer mesh.
 
\subsection{Influence of the Lam\'{e} parameters on the mesh quality}\label{subsection_InfluenceLame}
\begin{figure}
\centering
    \begin{subfigure}[b]{0.45\textwidth}
    \begin{overpic}[width = 0.9\textwidth]{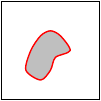}
    \put(38,40){$D_{\text{int}}$}
    \put(13,13){$D_0$}
    \put(44,72){\textcolor{red}{$\Gi $}}
    \put(-10,47){$\partial D$}
    \end{overpic}
\caption{Example of the domain $D$} 
\label{fig_D}
\end{subfigure}
~
\begin{subfigure}[b]{0.45\textwidth}
\includegraphics[scale=0.185]{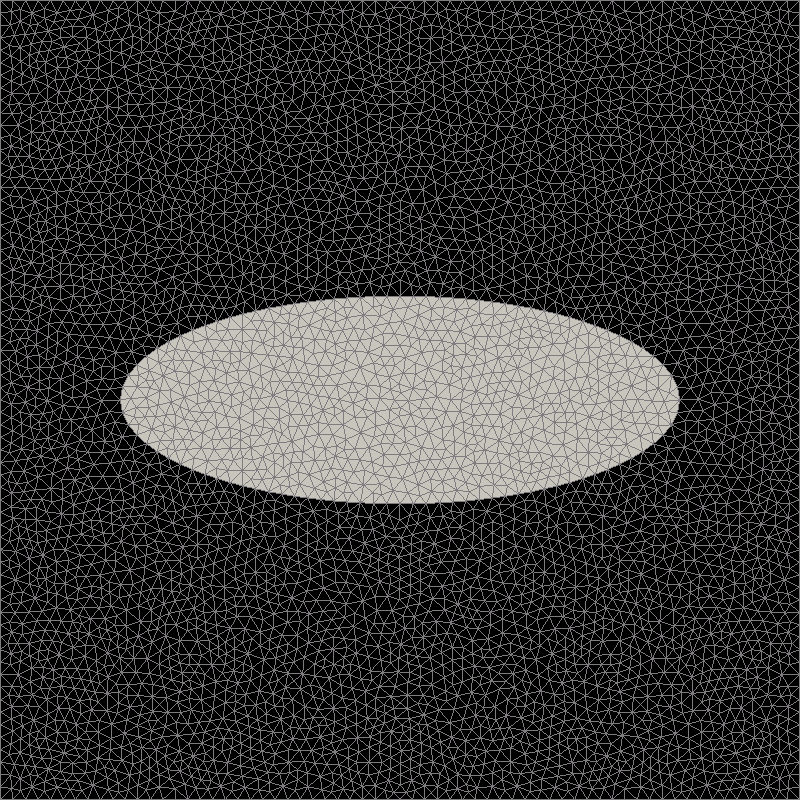}
\caption{Initial configuration $D^0$}
\label{fig:initial_experiment2}
\end{subfigure}
\caption{Experiment: Influence of Lam\'{e} parameters}
\end{figure}

For simplicity, from now on we only will consider shapes as depicted in~\Cref{fig_D}; i.e., consisting only of one interior interface. The main goal of this experiment is to show the importance of the parameters $\mu_{\min}$ and $\mu_{\max}$ in achieving better aspect ratios of the final meshes. The initial shape was chosen to be an ellipse as depicted in~\Cref{fig:initial_experiment2} on a mesh with about $10,000$ triangles. We used the Armijo line search rule~\eqref{eq:armijo-rule} with $\alpha = 400$, $\rho = 0.5$, $c = 10^{-4}$ and 
$\kappa_1\sim \mathcal{N}(1.5,10^{-2},1,2)$, $\kappa_{\text{int}}\sim \mathcal{N}(4,10^{-2},3,5)$ and $g \sim \mathcal{N}(10,10^{-2},9,11)$. We let the algorithm iterate 200 times.

\begin{figure}[ht]\centering
\begin{tikzpicture}[zoomboxarray,black and white=cycle]
    \node [image node]{\includegraphics[width=0.45\textwidth]		                      {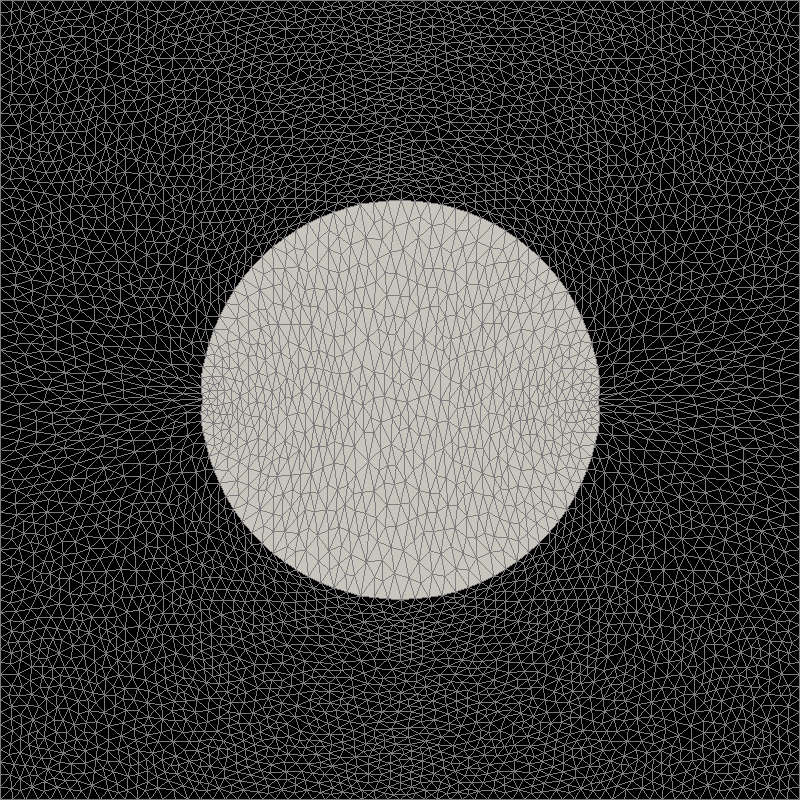} };
    \zoombox[magnification = 3]{0.04,0.5}
    \zoombox{0.5,0.75}
    \zoombox[magnification = 3]{0.5,0.04}
    \zoombox{0.75,0.5}
\end{tikzpicture}
\caption{(a) Shape obtained after 200 iterations using $\mu_{\min}=0.5$ and $\mu_{\max}=1$, (b)~zoomed areas}
\label{fig:mumin05_mumax1}
\end{figure}

\begin{figure}[ht]\centering
\begin{tikzpicture}[zoomboxarray,black and white=cycle]
    \node [image node]{\includegraphics[width=0.45\textwidth]		                      {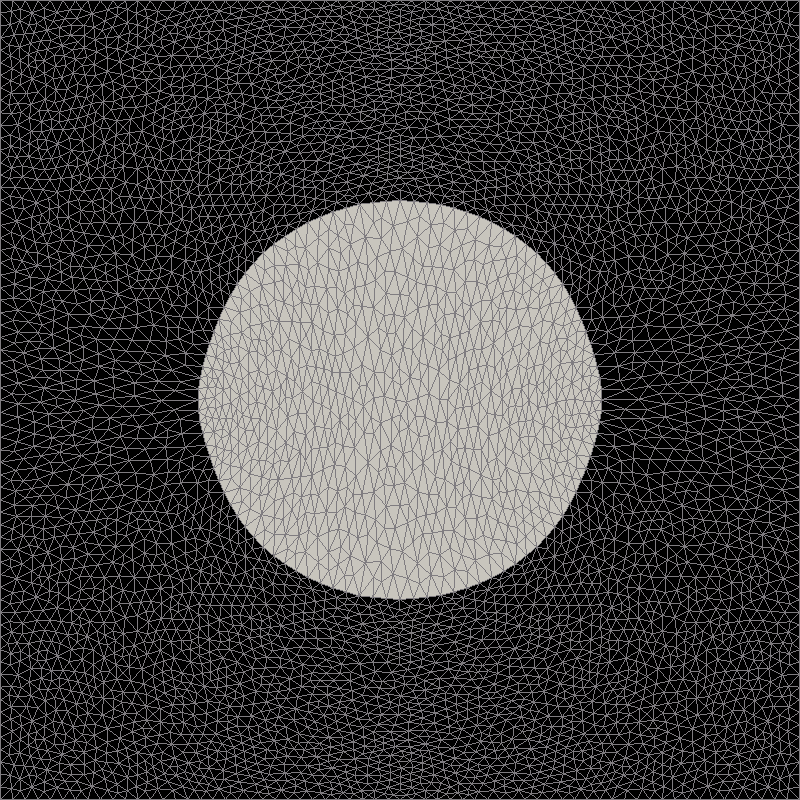} };
    \zoombox[magnification = 3]{0.04,0.5}
    \zoombox{0.5,0.75}
    \zoombox[magnification = 3]{0.5,0.04}
    \zoombox{0.75,0.5}
\end{tikzpicture}
\caption{(a) Shape obtained after 200 iterations using $\mu_{\min}=10$ and $\mu_{\max}=25$, (b)~zoomed areas}
\label{fig:mumin10_mumax25}
\end{figure}

\Cref{fig:mumin05_mumax1,fig:mumin10_mumax25} show that of the tested values, the aspect ratio for the final mesh is the best in the case where we choose $\mu_{\min} = 10, \mu_{\max}=25$. These results suggest that the correct choice of these values is problem-dependent and would be an opportunity for future research. Specifically, we have improved the aspect ratios of the cells near to the moving boundary by increasing the value of $\mu_{\min}$ and the  non-moving boundary by increasing the value of $\mu_{\max}$.   

\begin{figure}[ht]\centering
\begin{tikzpicture}[zoomboxarray,black and white=cycle]
    \node [image node]{\includegraphics[width=0.45\textwidth]		                      {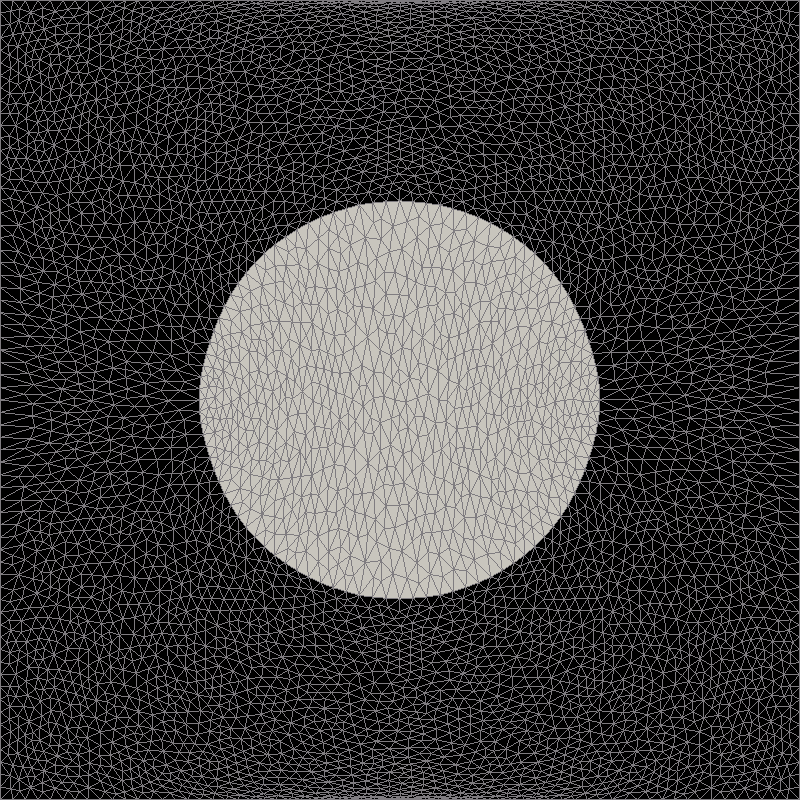} };
    \zoombox[magnification = 3]{0.04,0.5}
    \zoombox{0.5,0.75}
    \zoombox[magnification = 3]{0.5,0.04}
    \zoombox{0.75,0.5}
\end{tikzpicture}
\caption{(a) Shape obtained after 200 iterations using $\mu_{\min}=1$ and $\mu_{\max}=25$, (b)~zoomed areas}
\label{fig:mumin1_mumax25}
\end{figure}  

\subsection{Different target shapes}
\begin{figure}
  \begin{subfigure}[b]{0.33\linewidth}
    \includegraphics[height = 3.1cm,  width = 3.1cm]{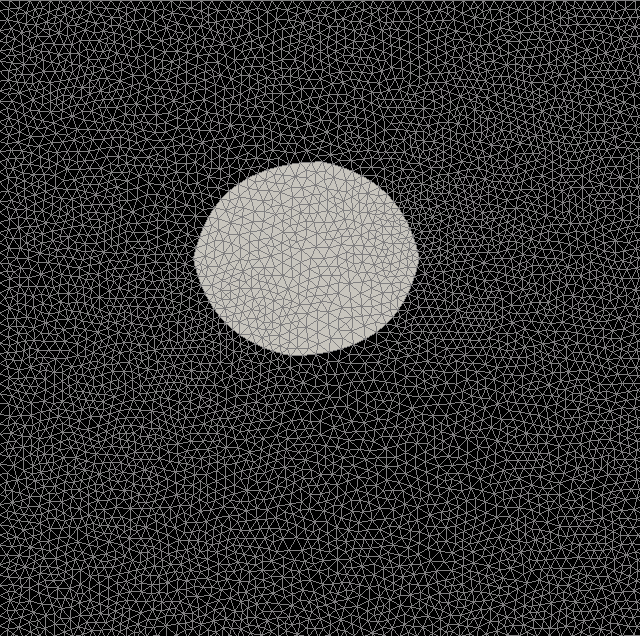}
    \caption{Initial configuration $D^0$}  
    \label{fig:circle-to-ellipse_target}
  \end{subfigure}%
  \begin{subfigure}[b]{0.33\linewidth}
    \centering
    \includegraphics[height = 3.1cm,  width = 3.1cm]{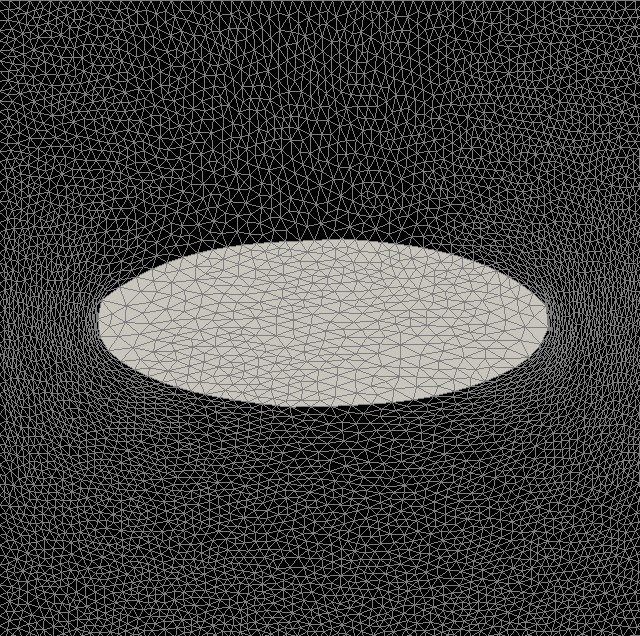}   
    \caption{$D^{200}$}
  \end{subfigure}%
  \begin{subfigure}[b]{0.33\linewidth}
    \centering
    \includegraphics[height = 3.1cm]{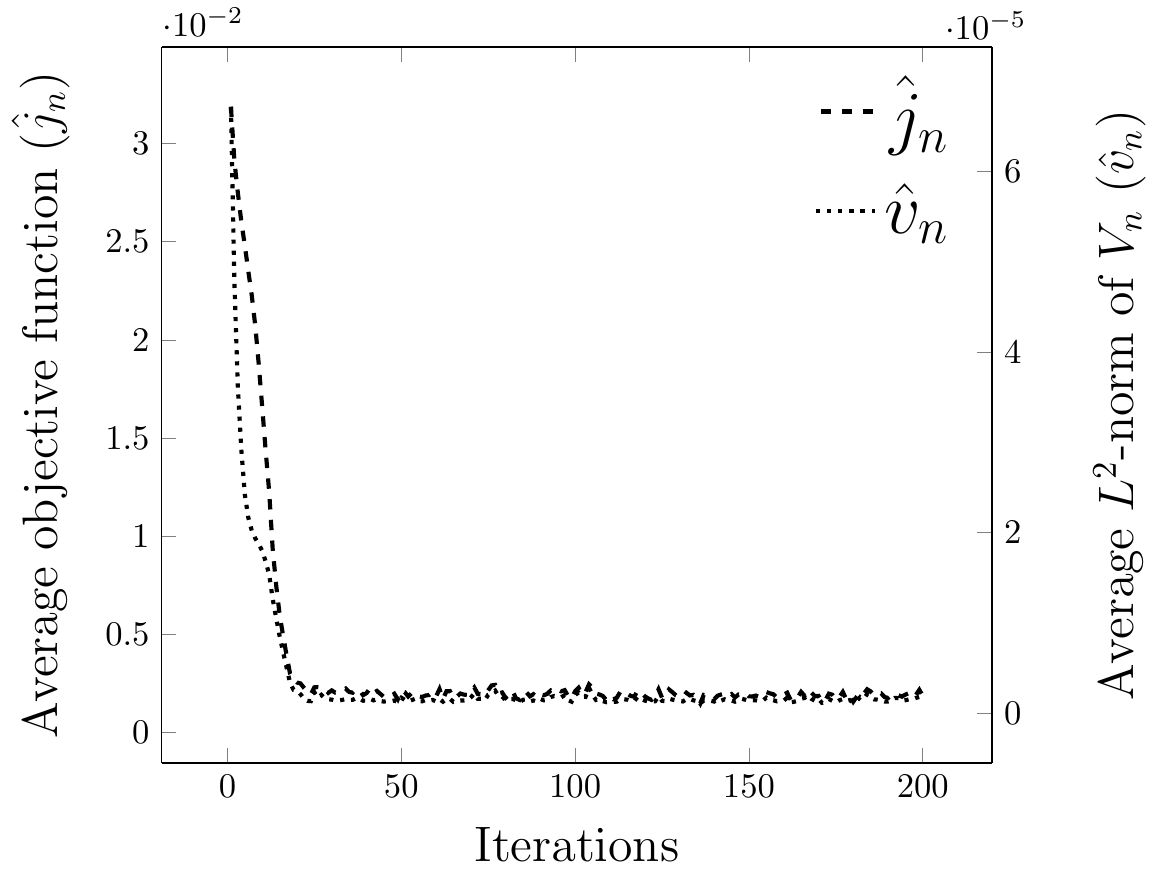}   		
    \caption{Convergence plot}
  \end{subfigure}%
  \caption{Convergence of circle to ellipse}
  \label{fig:circle-to-ellipse}
\end{figure}

\begin{figure}
  \begin{subfigure}[b]{0.33\linewidth}
    \includegraphics[height = 3.1cm,  width = 3.1cm]{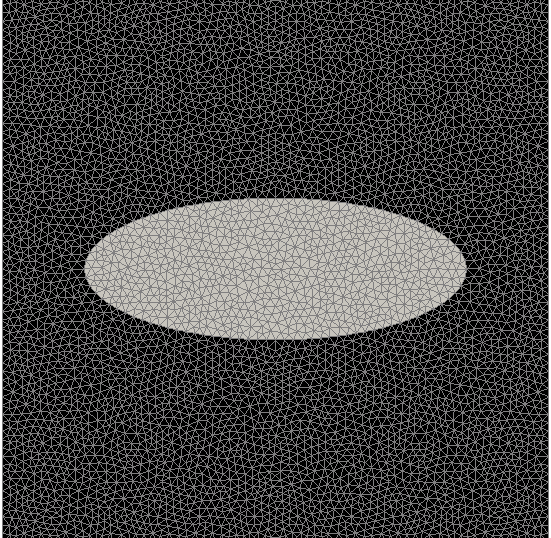}
    \caption{Initial configuration $D^0$}  
    \label{fig:ellipse-to-circle_target}
  \end{subfigure}%
  \begin{subfigure}[b]{0.33\linewidth}
    \centering
    \includegraphics[height = 3.1cm,  width = 3.1cm]{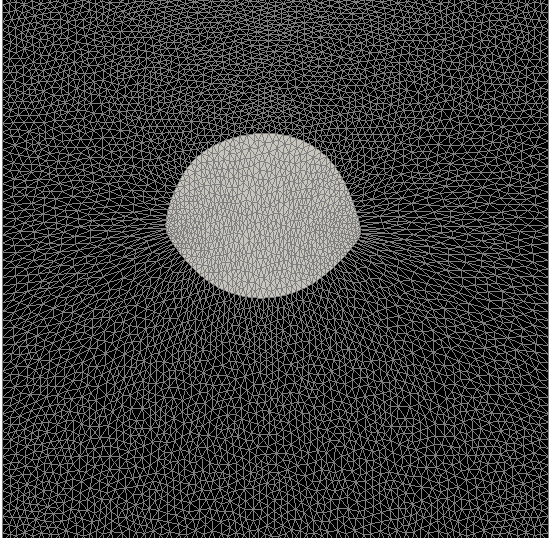}   
    \caption{$D^{200}$}
  \end{subfigure}%
  \begin{subfigure}[b]{0.33\linewidth}
    \centering
    \includegraphics[height = 3.1cm]{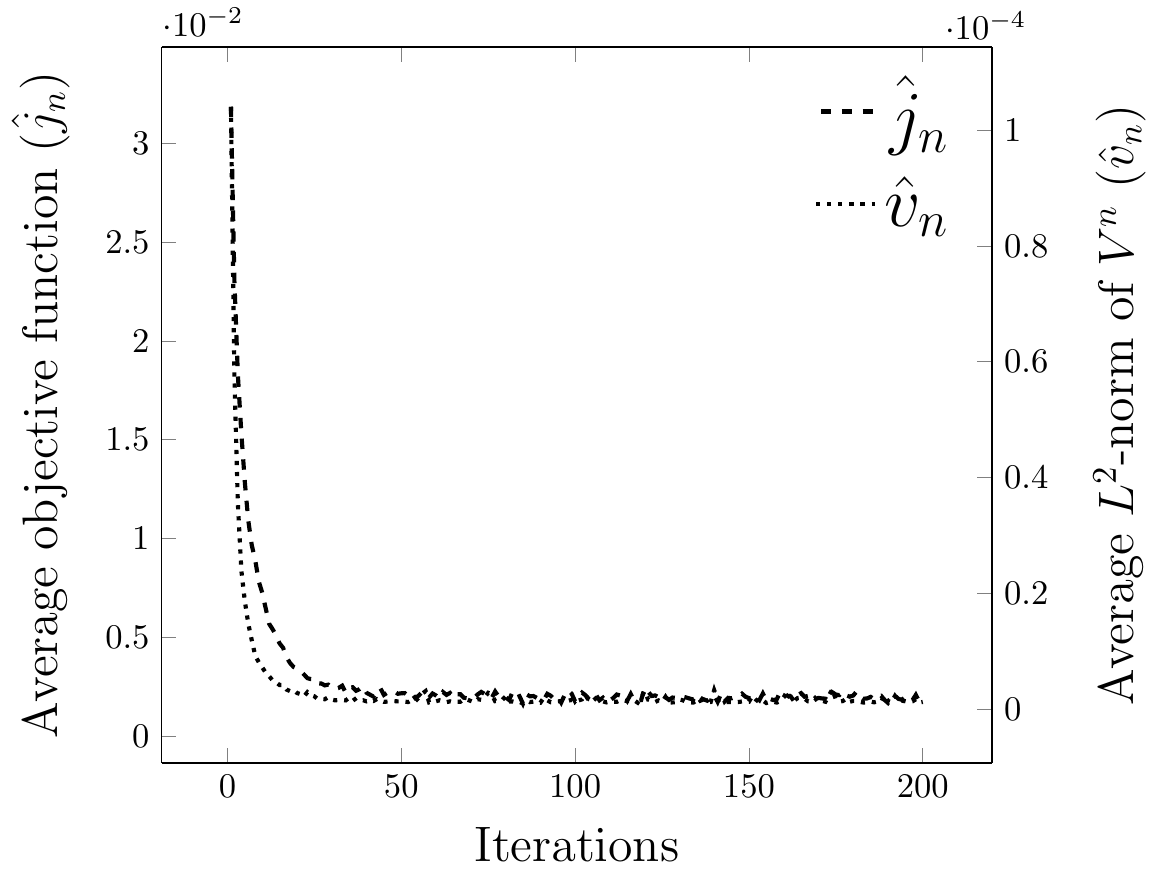}   		
    \caption{Convergence plot}
  \end{subfigure}%
  \caption{Convergence of ellipse to circle}
  \label{fig:ellipse-to-circle}
\end{figure}

In this experiment, we show convergence to different target shapes.  We used the Armijo line search rule~\eqref{eq:armijo-rule} with $\alpha = 400$, $\rho = 0.5$, $c = 10^{-4}$ and $\kappa_0  \sim \mathcal{N}(1.5,10^{-2},1,2) $, $ \kappa_{\text{int}}\sim \mathcal{N}(4,10^{-2},3,5)$ and $g \sim \mathcal{N}(10,10^{-2},9,11)$.
The measurement $\bar{y}$ was generated for two different target shapes $\bar{u}$: the ellipse (\Cref{fig:ellipse-to-circle_target}) and a small circle (\Cref{fig:circle-to-ellipse_target}). The starting shape $D^0$, final shape $D^{200}$ and convergence plot is shown in each figure. In the convergence plots, one sees a clear decrease in the objective function value and the norm of the shape gradient as well as significant progress toward to the target shape.

\subsection{Influence of the individual random variables}
\label{sec:individual_variables}

\begin{figure}
\centering
\begin{subfigure}[b]{0.45\textwidth}
\includegraphics[scale=0.27]{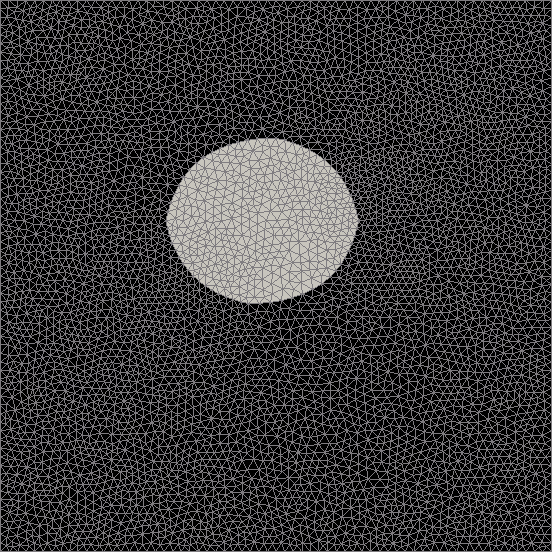}
\caption{Initial configuration $D^0$} 
\label{fig:initial_experiment4}
\end{subfigure}
~
\begin{subfigure}[b]{0.45\textwidth}
\includegraphics[scale=0.27]{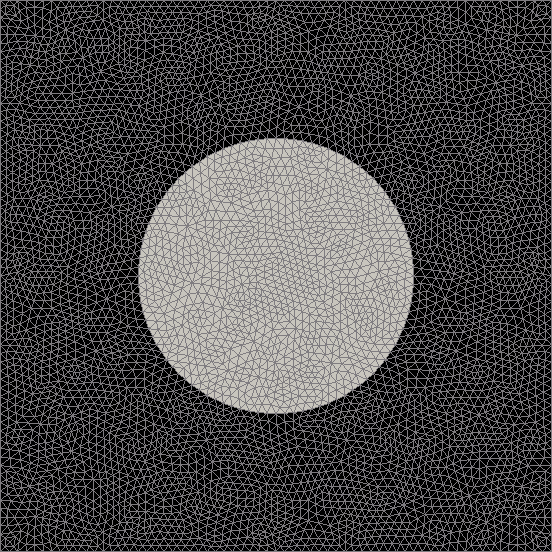}
\caption{Target shape $\bar{u}$} 
\label{fig:final_experiment4}
\end{subfigure}
\caption{Experiment: Influence of the individual random variables} 
\end{figure}

The main purpose of this experiment is to study the influence of the individual random variables on the performance of the algorithm. That means we will choose different values to be deterministic while other will remain random.
The initial shape $D^0$ is depicted in~\Cref{fig:initial_experiment4} and the target shape $\bar{u}$ is depicted in~\Cref{fig:final_experiment4}.

\begin{figure}
\centering
\begin{subfigure}[b]{0.45\textwidth}
\includegraphics[scale=0.5]{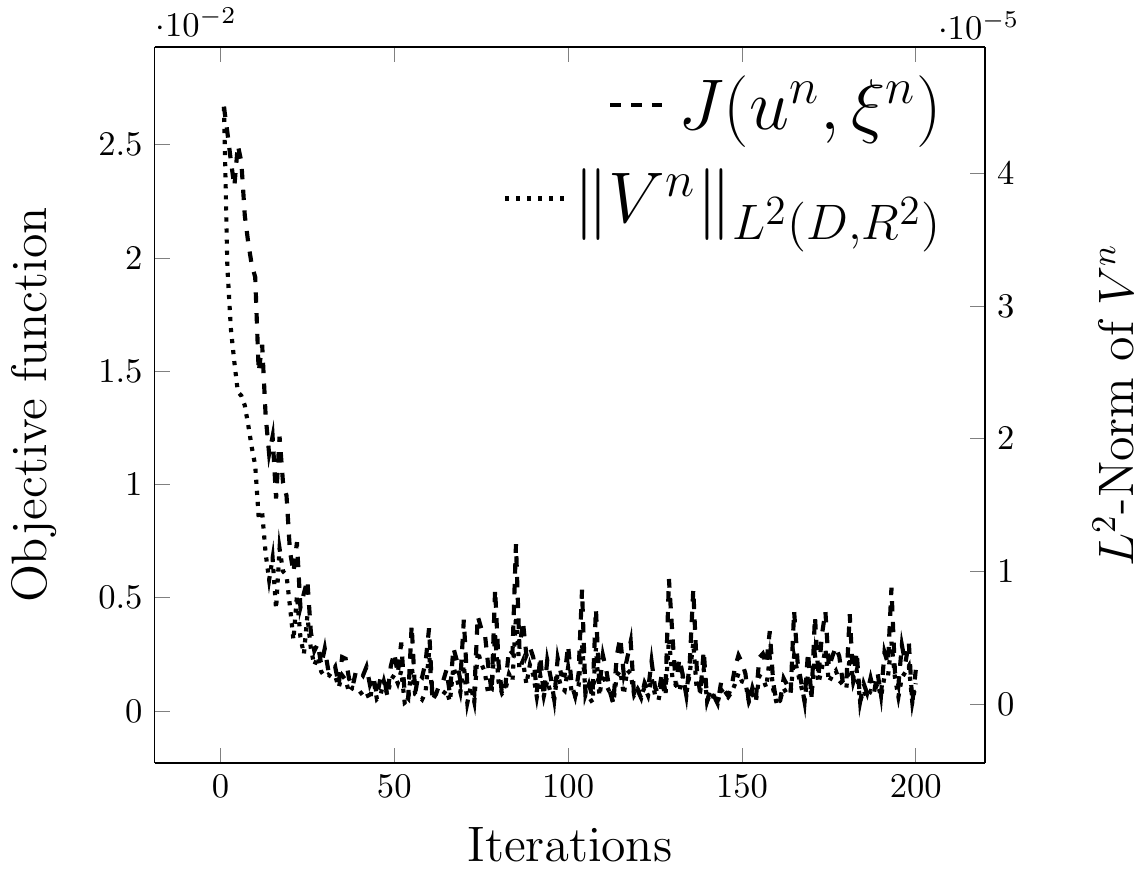}
\caption{ $\kappa_0$ random} 
\label{fig:k1_random}
\end{subfigure}
~
\begin{subfigure}[b]{0.45\textwidth}
\includegraphics[scale=0.5]{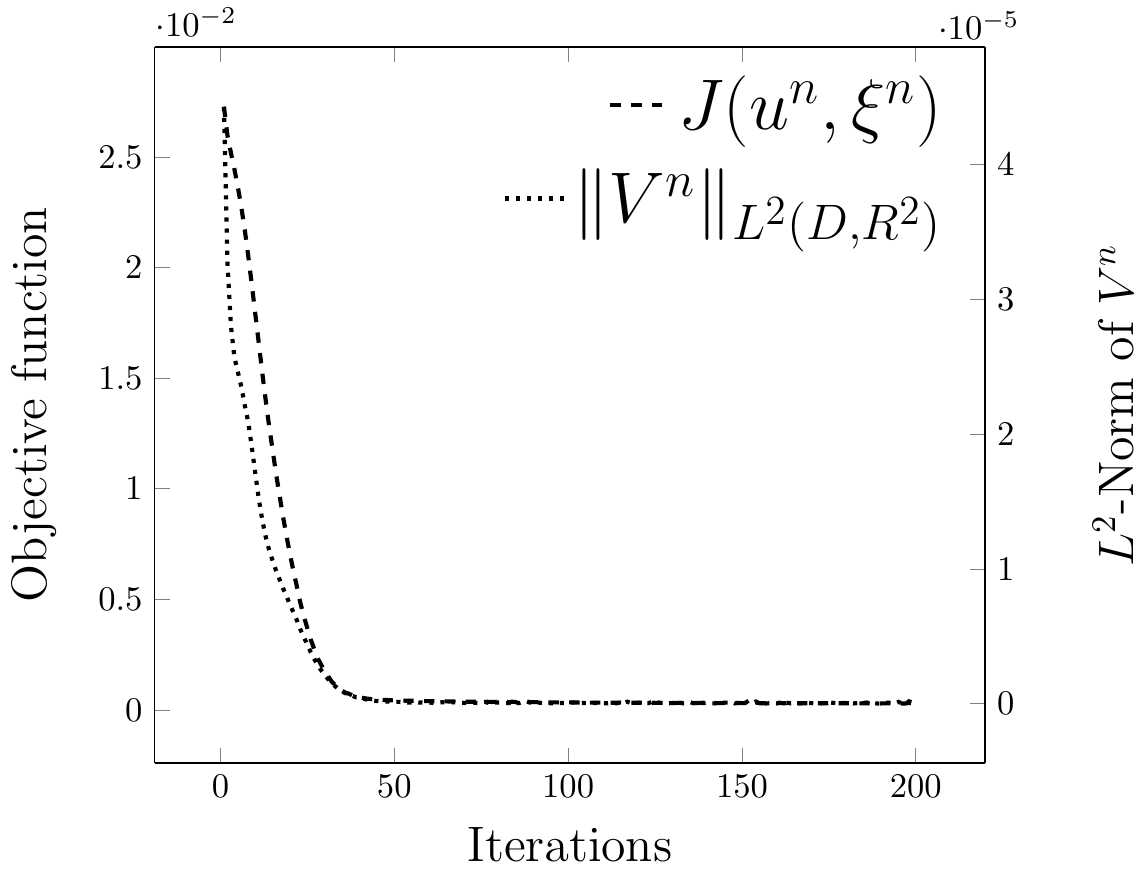}
\caption{$\kappa_{\text{int}}$ random} 
\label{fig:k2_random}
\end{subfigure}
~
\begin{subfigure}[b]{0.45\textwidth}
\includegraphics[scale=0.5]{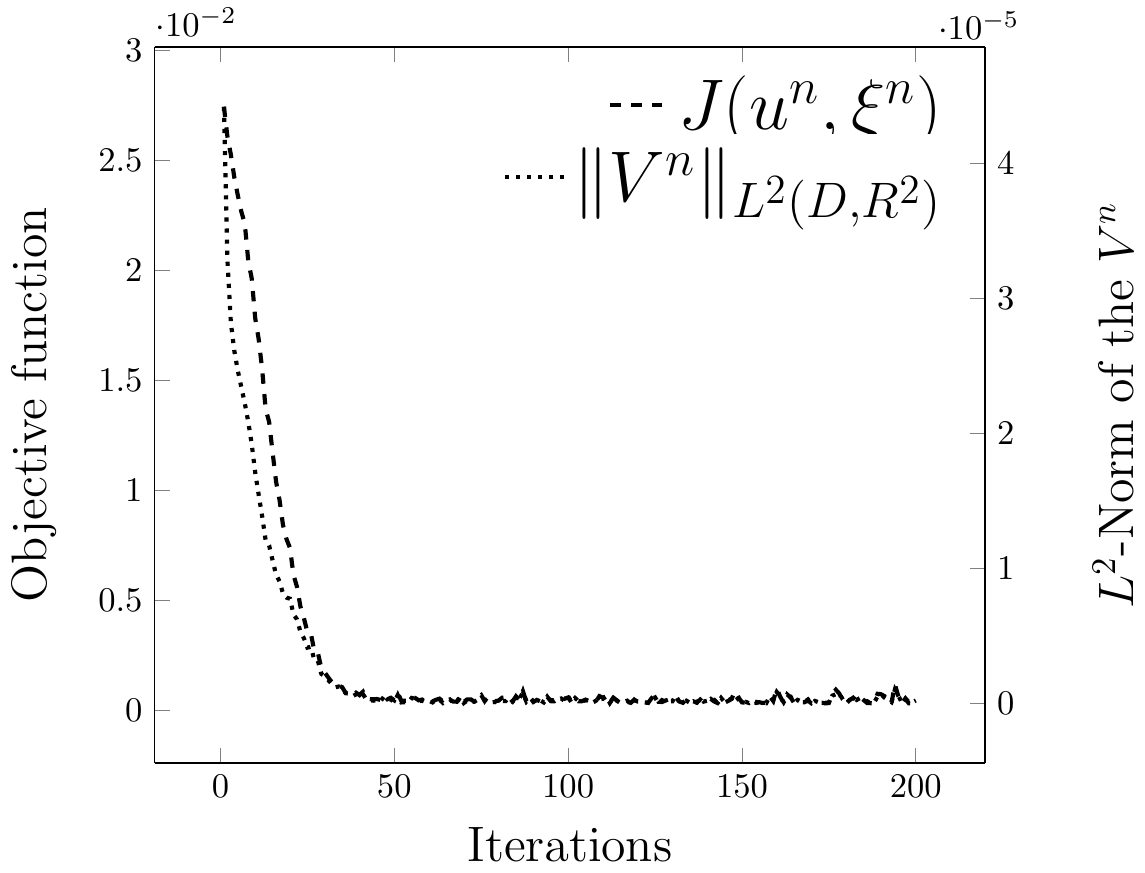}
\caption{ $g$ is random} 
\label{fig:g_random}
\end{subfigure}
~
\begin{subfigure}[b]{0.45\textwidth}
\includegraphics[scale=0.5]{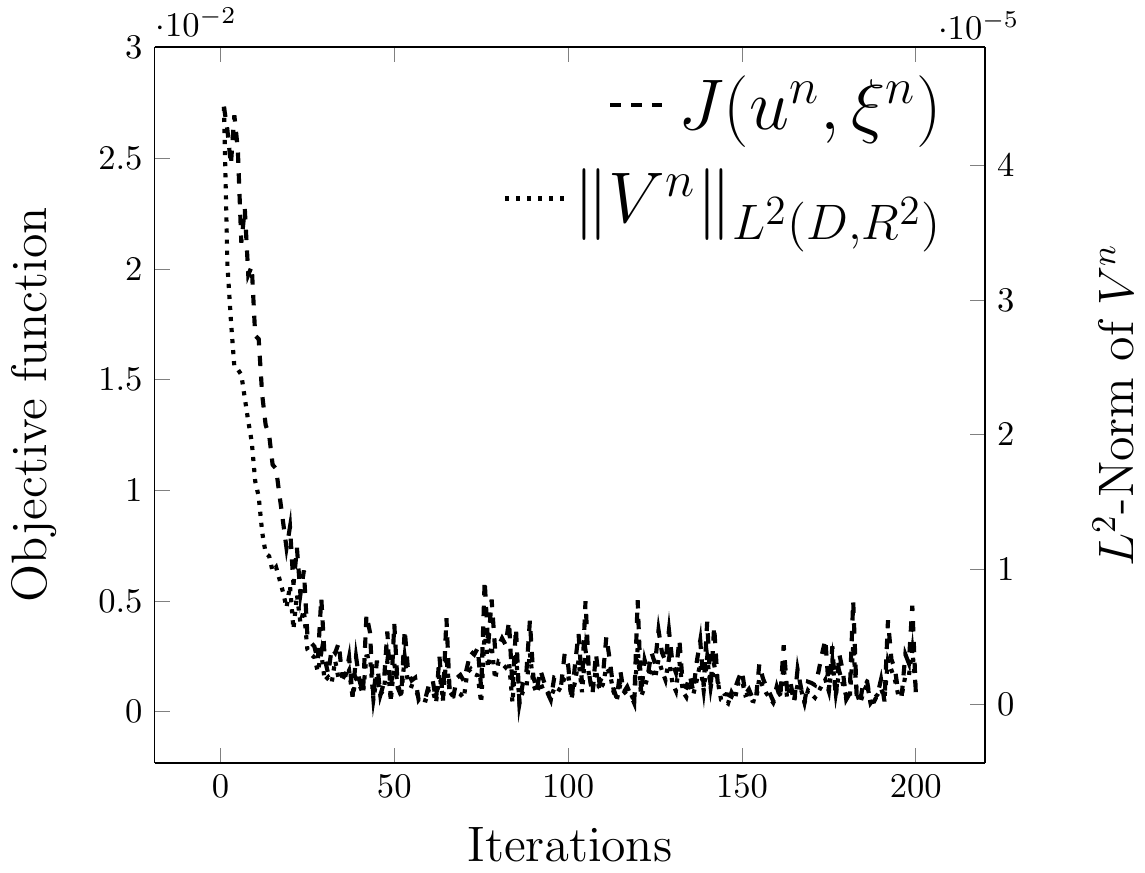}
\caption{$\kappa_0$ and $\kappa_{\text{int}}$ random}
\label{fig:k1_and_k2_random}
\end{subfigure}
~
\begin{subfigure}[b]{0.45\textwidth}
\includegraphics[scale=0.5]{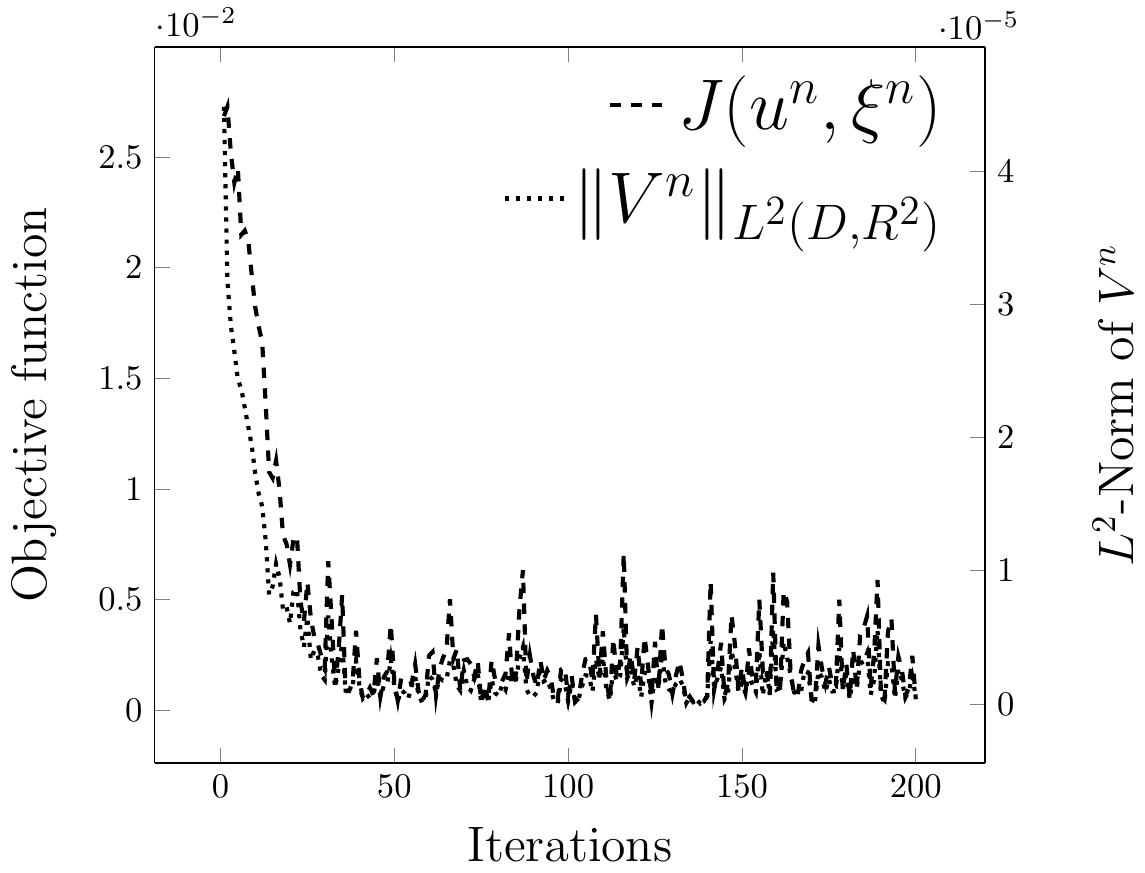}
\caption{Completely random}
\label{fig:random}
\end{subfigure}
~
\begin{subfigure}[b]{0.45\textwidth}
\includegraphics[scale=0.5]{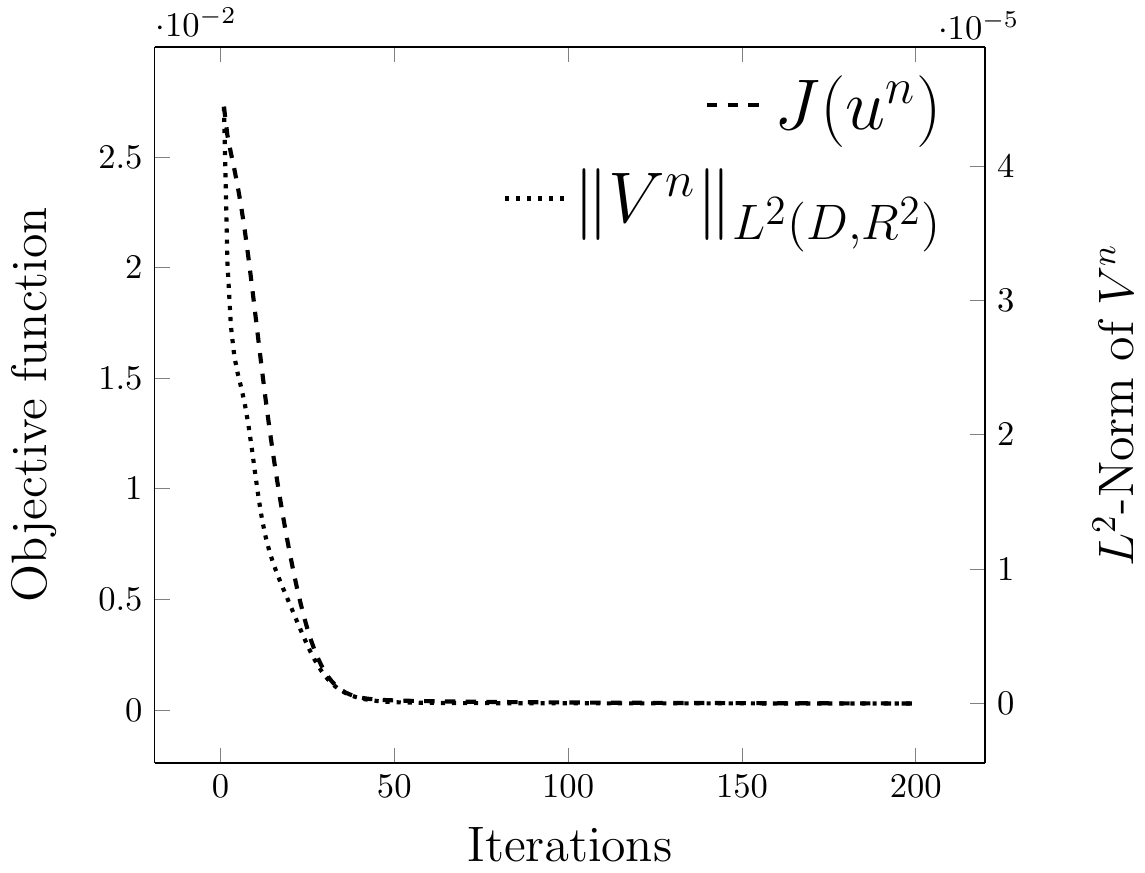}
\caption{ Deterministic } 
\label{fig:deterministic}
\end{subfigure}

\caption{Convergence plots}
\label{fig:experiment4}
\end{figure}

The behavior of the objective function and the norm of the shape gradient are depicted in~\Cref{fig:experiment4}, where we choose the following parameters. 
\begin{description}
\item[Figure~\ref{fig:k1_random}:] $\kappa_0\sim \mathcal{N}(1.5,10^{-2},1,2)$, while $\kappa_{\text{int}} = 4$ and $g = 10$.
\item[Figure~\ref{fig:k2_random}:] $\kappa_{\text{int}}\sim \mathcal{N}(4,10^{-2},3,5)$, while $\kappa_0 = 1.5$ and $g = 10$.
\item[Figure~\ref{fig:g_random}:]$g \sim \mathcal{N}(10,10^{-2},9,11)$, while $\kappa_0 = 1.5$ and $\kappa_{\text{int}} = 4$.
\item[Figure~\ref{fig:k1_and_k2_random}:] $\kappa_0\sim \mathcal{N}(1.5,10^{-2},1,2)$, $\kappa_{\text{int}}\sim \mathcal{N}(4,10^{-2},3,5)$, and $g = 10$. 
\item[Figure~\ref{fig:random}:] $\kappa_0\sim \mathcal{N}(1.5,10^{-2},1,2)$, $\kappa_{\text{int}}\sim \mathcal{N}(4,10^{-2},3,5)$ and $g \sim \mathcal{N}(10,10^{-2},9,11)$.
\item[Figure~\ref{fig:deterministic}:] $\kappa_0 = 1.5$, $\kappa_{\text{int}}= 4 $ and $g=10$. 
\end{description}

By comparing~\Cref{fig:k2_random,fig:deterministic} we can conclude that the randomness in $\kappa_{\text{int}}$ influences the behavior of the algorithm only minimally. That agrees with the model, since this parameter affects only the moving boundary from the inside. 
Comparing~\Cref{fig:k1_random,fig:k1_and_k2_random,fig:random}, we observe that the stochastic shape gradient appears to be most sensitive to the parameter $\kappa_0$, which induces greater oscillations in the objective function throughout optimization. We can also see in~\Cref{fig:g_random} that the stochastic shape gradient and objective function value do not appear to be very sensitive to randomness in the boundary condition $g$.

\subsection{Step size rules for high variance distributions}

\begin{figure}
\centering
\begin{subfigure}[b]{0.45\textwidth}
\centering
\includegraphics[height=4.3cm]{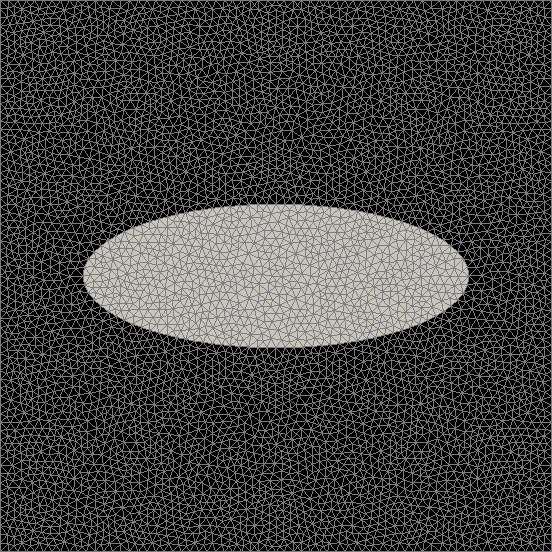}
\caption{Initial configuration $D^0$} 
\label{fig:initial_experiment5}
\end{subfigure}
~
\begin{subfigure}[b]{0.45\textwidth}
\centering
\includegraphics[height=4.3cm]{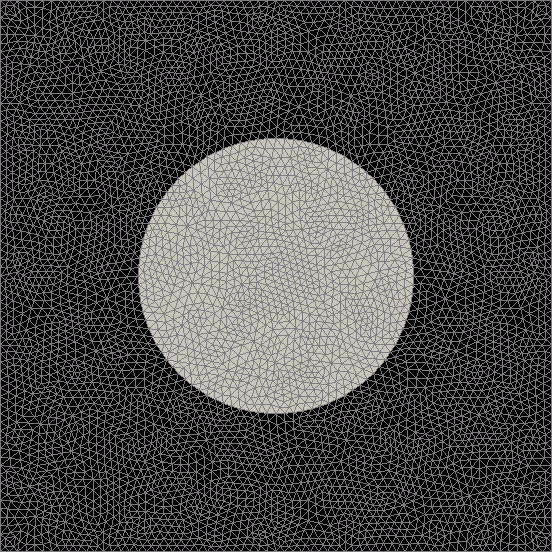}
\caption{Target shape $\bar{u}$} 
\label{fig:final_experiment5}
\end{subfigure}
\caption{Experiment: Step size rules} 
\end{figure}

In this experiment, we will show the influence of the step size $t^n$ in~\Cref{alg:stochastic_gradient_shape_optimization} on convergence in the high-variance case. We will compare the Robbins-Monro step size rule~\eqref{eq:Robbins-Monro} with the Armijo line search rule~\eqref{eq:armijo-rule} as well as a modification of the latter. To demonstrate the effect of high variance on the problem, we will work with the following random variables with higher variance than in the past examples: $\kappa_0  \sim \mathcal{N}(1.5,0.2,1,2) $, $ \kappa_{\text{int}}\sim \mathcal{N}(4,0.2,3,5)$ and $g \sim \mathcal{N}(10,0.2,9,11)$. The initial shape is shown in~\Cref{fig:initial_experiment5} and the target measurement $\bar{y}$ is generated on the shape in~\Cref{fig:final_experiment5}. Convergence plots are shown with $m=1000$ samples for each iteration $n$.

\subsubsection{Robbins-Monro step size}
\begin{figure}
  \begin{subfigure}[b]{0.45\linewidth}
  \centering
    \includegraphics[height=4.3cm]{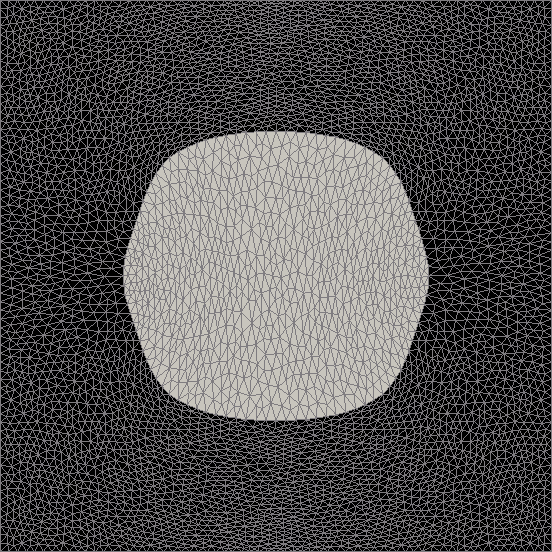}
    \caption{$D^{200}$}  
  \end{subfigure}%
  \begin{subfigure}[b]{0.45\linewidth}
    \centering
    \includegraphics[height=4.3cm]{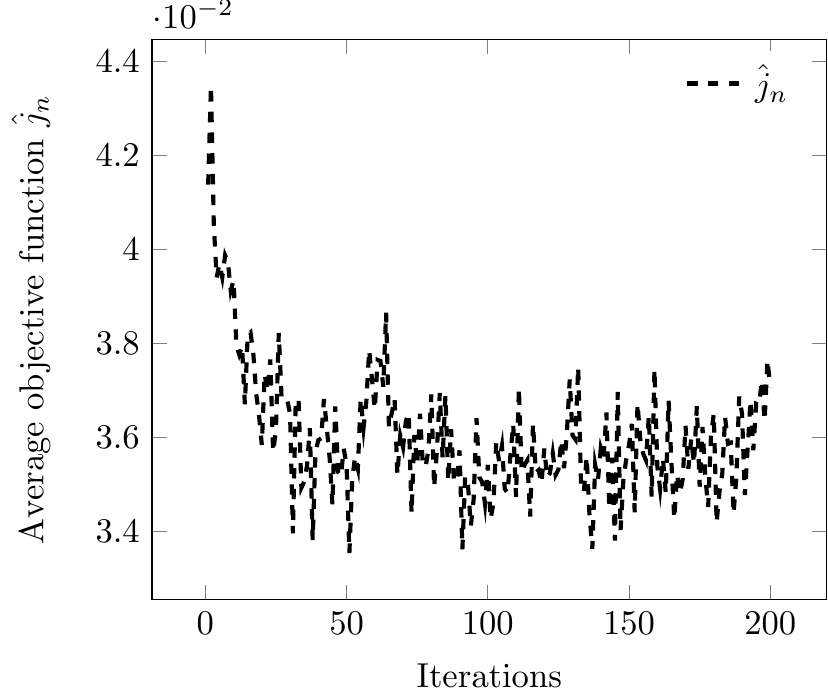}   
    \caption{Objective function}
  \end{subfigure}%
  \caption{Robbins-Monro rule $t^n = 800 n^{-0.85}$}
  \label{fig:RM-800}
\end{figure}

\begin{figure}
  \begin{subfigure}[b]{0.45\linewidth}
  \centering
    \includegraphics[height=4.3cm]{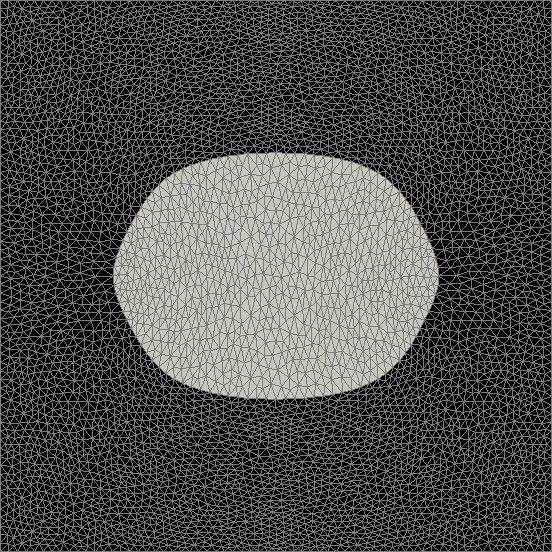}
    \caption{$D^{200}$}  
  \end{subfigure}%
  \begin{subfigure}[b]{0.45\linewidth}
    \centering
    \includegraphics[height=4.3cm]{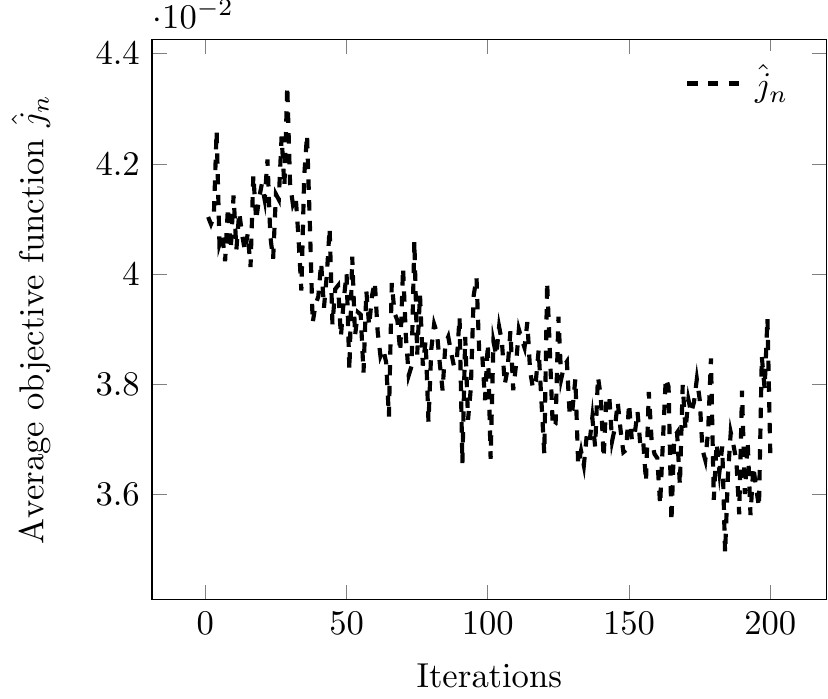}   
    \caption{Objective function}
  \end{subfigure}%
  \caption{Robbins-Monro rule $t^n = 400 n^{-0.85}$}
  \label{fig:RM-400}
\end{figure}

We work with the rule $t^n = \alpha n^{-0.85}$, which obviously satisfies the conditions in~\eqref{eq:Robbins-Monro}. In~\Cref{fig:RM-800,fig:RM-400}, we compare results of a simulation using $\alpha=800$ and $\alpha=400$, respectively. This demonstrates that a poorly chosen constant can drastically affect convergence speed. Interestingly, at $n=200$, the better choice $\alpha=800$ has $\hat{j}_{200} \approx 0.371,$ while $\alpha=400$ has $\hat{j}_{200} \approx 0.368$, so a slightly lower objective function value. However, as one clearly sees from the plots, even with $m=1000$ samples, the estimated function values are quite noisy; therefore, we cannot conclude that the shape obtained in~\Cref{fig:RM-400} has a lower objective function value than the one in~\Cref{fig:RM-800}. 

\subsubsection{Armijo line search}
\begin{figure}
  \begin{subfigure}[b]{0.45\linewidth}
  \centering
    \includegraphics[height = 4.3cm]{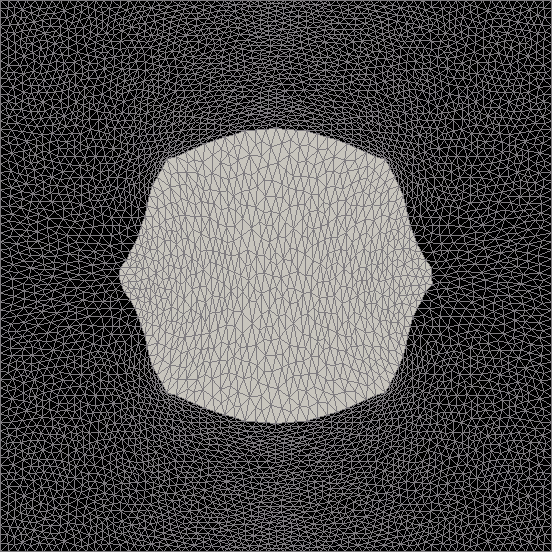}
    \caption{$D^{200}$}  
  \end{subfigure}%
  \begin{subfigure}[b]{0.45\linewidth}
    \centering
    \includegraphics[height=4.3cm]{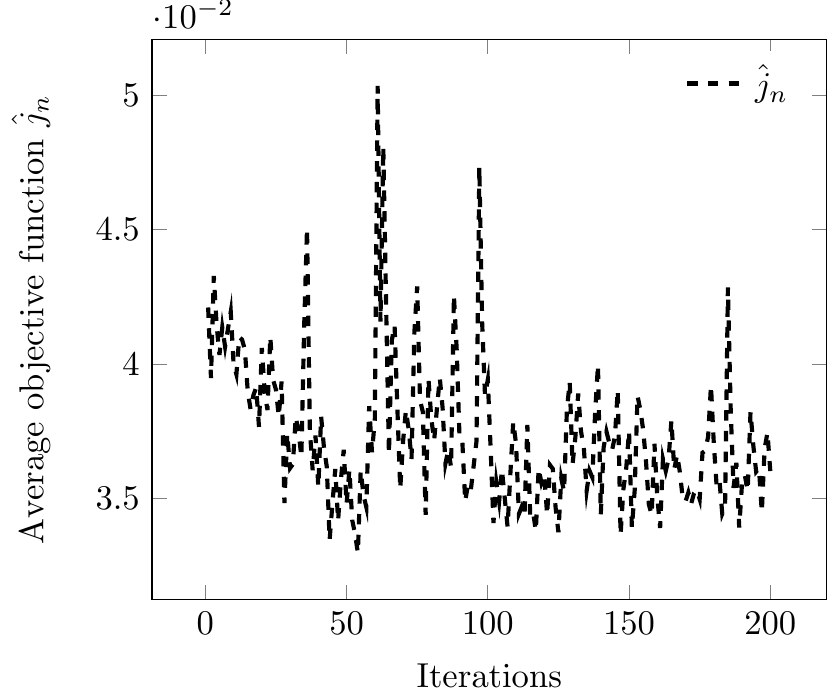}   
    \caption{Objective function}
  \end{subfigure}%
  \caption{Armijo line search with $t^n = 300 \cdot 0.5^{m_n}$}
  \label{fig:Armijo-300}
\end{figure}

We use the Armijo line search rule~\eqref{eq:armijo-rule} with $\alpha = 300$, $\rho = 0.5$, $c = 10^{-4}$ and $N_n \equiv 1$ sample per iteration. In~\Cref{fig:Armijo-300}, one sees the final shape as well as a plot of the objective function. In this simulation, the Armijo line search was successful in the first iteration for every $n$, effectively making the step size constant $t^n \equiv 150$. Unlike the Robbins-Monro rule, this step size rule does not reduce variance well, since large steps are taken even if one has already reached the neighborhood of the optimum. This induces large oscillations in the shape, which may compromise mesh integrity. The objective function value also oscillates heavily as a result, as one can see in the plot of the objective function. At $n=200$, $\hat{j}_{200} \approx 0.360$, where it is noted that one would need to take more samples before making a fair comparison with the previous example.

\subsubsection{Armijo line search with damping}
\begin{figure}
  \begin{subfigure}[b]{0.45\linewidth}
  \centering
    \includegraphics[height = 4cm]{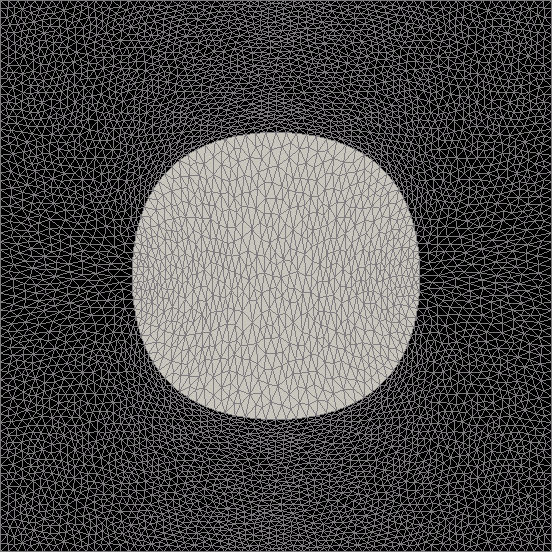}
    \caption{$D^{200}$}  
  \end{subfigure}%
  \begin{subfigure}[b]{0.45\linewidth}
    \centering
    \includegraphics[height = 4cm]{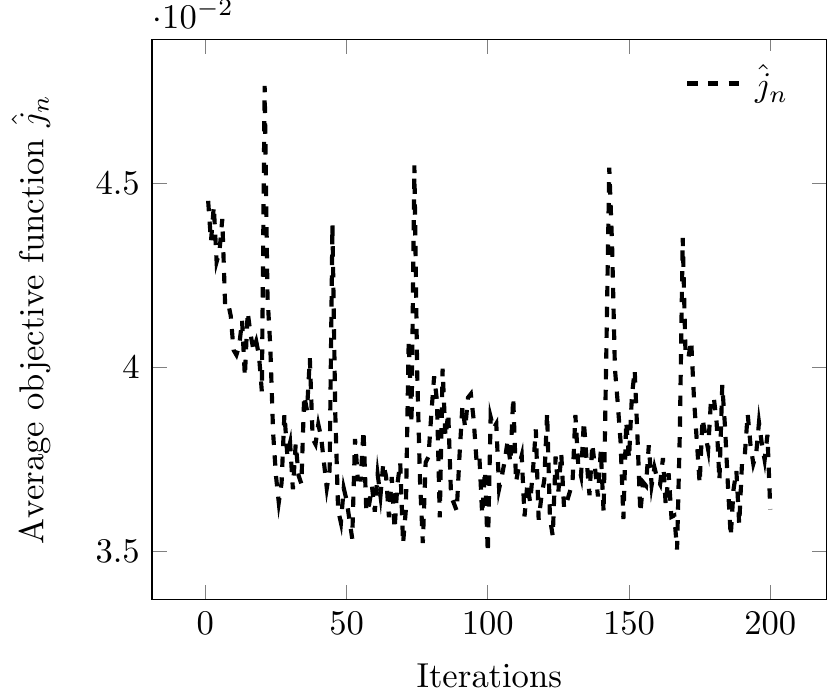}   
    \caption{Objective function}
  \end{subfigure}%
  \caption{Armijo line search with $t^n = \alpha_n \cdot 0.5^{m_n}$}
  \label{fig:convergence-2e-1}
\end{figure}

To mimic the variance reduction property of the Robbins-Monro rule, we introduce Armijo line search with damping
to reduce oscillations around the optimum. This line search takes the form
\begin{equation}
 \label{eq:linesearch-damping}
 t^n = \alpha_n \rho^{m_n}
\end{equation}
for a decreasing sequence $\alpha_n \in [\alpha_{\min}, \alpha_{\max}].$ Here, we choose $\alpha_0 = 400,$ $\rho = 0.5$, and $c = 10^{-4}$, and damp $\alpha_n$ every $20^{th}$ iteration (i.e. $\alpha_{20} = 0.9\alpha_0$, $\alpha_{21} =\cdots = \alpha_{39} = \alpha_{20}$, $\alpha_{40} = 0.9\alpha_{20}$, etc.) \Cref{fig:convergence-2e-1} shows the shape and a plot of objective function. Although the function values still oscillate, the extreme values are lower. At $n=200$, $\hat{j}_{200} \approx 0.361$. While this experiment shows some promise in the hybrid step size rule~\eqref{eq:linesearch-damping}, it is to be noted that to the authors' knowledge, no convergence theory for this method exists; this would be an opportunity for future research. 

\subsubsection{Concluding observations}

\begin{figure}[ht]\centering
\begin{tikzpicture}[zoomboxarray, zoomboxarray columns=2,
    zoomboxarray rows=1,black and white=cycle]
    \node [image node]{\includegraphics[width=0.35\textwidth]		                      {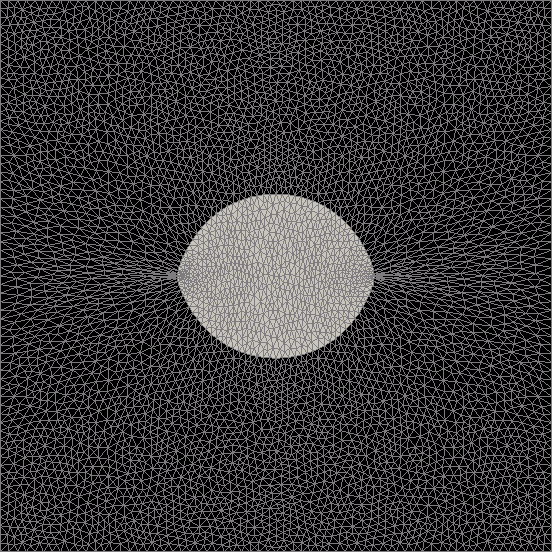} };
    \zoombox{0.3,0.5}
   	\zoombox{0.7,0.5} 
\end{tikzpicture}
\caption{(a) $D^{2}$ with Robbins-Monro rule~\eqref{eq:Robbins-Monro} and $\alpha=2000$, (b)~zoomed areas}
\label{fig:RM-broken-mesh}
\end{figure}

\begin{figure}[ht]\centering
\begin{tikzpicture}[zoomboxarray,black and white=cycle]
    \node [image node]{\includegraphics[width=0.35\textwidth]		                      {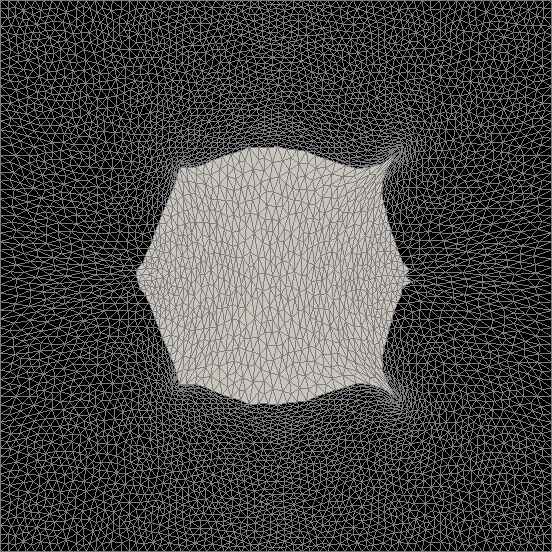} };
    \zoombox{0.35,0.7}
   	\zoombox{0.7,0.3} 
   	\zoombox{0.35,0.3}
    \zoombox{0.7,0.7}
\end{tikzpicture}
\caption{(a) $D^{200}$ with Armijo line search rule \eqref{eq:armijo-rule} and $\alpha=400$, (b)~zoomed areas}
\label{fig:Armijo-broken-mesh}
\end{figure}

Notably, all of these methods are quite sensitive to the scaling of the step size, requiring the user to tune parameters offline. This is in some sense unavoidable for iterative methods on nonconvex problems. For all step size rules, if the scaling is chosen to be too large, mesh destruction or bad aspect ratios of the mesh may appear as seen in~\Cref{fig:RM-broken-mesh,fig:Armijo-broken-mesh}. For the Robbins-Monro rule, the danger of mesh destruction occurs at the beginning, when the step size is its largest. \Cref{fig:RM-broken-mesh} shows how the mesh has bad aspect ratios already at the second iteration if $t^n = 2000\cdot{n}^{-0.85}$. \Cref{fig:Armijo-broken-mesh} shows a mesh with bad aspect ratios at the final iteration, where high oscillations around iteration $n=175$ caused a compromised mesh.

The Armijo rule, which is widely used for deterministic problems, should be used with care for stochastic problems. As mentioned in Section~\ref{sec:algorithmic_details}, current convergence theory requires that one takes more and more samples of the gradient to achieve variance reduction. The Robbins-Monro rule, which naturally reduces variance over the optimization, suffers from a lack of robustness: slight changes in the scaling of the step size lead to drastically different efficiency in the first stage of optimization, as we saw in the experiments. An idea is to try to reduce variance while still allowing for large step sizes, which we used to develop the damping procedure. This method, however, introduces an additional parameter needing to be tuned, namely the damping rate. So far, no theoretical results exist for such a procedure, but it would be an interesting topic to explore, particularly for nonconvex problems.

\section{Conclusion}
\label{sec:conclusion}
A stochastic shape interface problem was investigated computationally. 
The problem under consideration is of high importance because adding randomness makes numerical simulations for real-world problems with uncertainty in the model.
In order to solve the problem, we combine the gradient descent method using the shape derivative in its volume expression with the well-known stochastic gradient, which results in a stochastic gradient descent method for shape optimization problems in shape spaces. We present various numerical experiments in this paper by applying this method to the shape interface model.
Numerical solutions of the problem are presented, where the domain consists of multiple shapes.
Moreover, the influence of the Lam\'{e} parameters on the mesh quality and the influence of individual random variables on the performance of the algorithm are investigated.
In a further experiment, we perform computations for different target shapes, observing the convergence of a circle to an ellipse and vice versa. In a final experiment, which focuses on the various step size rules, we use an example with higher variance and compare simulations for the Robbins-Monro step size rule, the Armijo line search rule, and a hybrid method we call Armijo line search with damping. The strengths and weaknesses of the methods are compared.

Many interesting questions remain, some of which we will attempt to address in future work. 
For example, how to extend our method efficiently to problems with two-dimensional shapes, i.e., where the underlying domains are subsets of $\R^3$. Moreover, it is well done that with highly dimensional problems it arises the question of improving the performance of the algorithm. For example, parallelization techniques could be solution.
Additionally, it would be useful to compute efficiency estimates so that a termination condition can be defined for the algorithm. 
As we demonstrated in the last experiment, step size plays a big role, especially when the variance is high;  although no convergence theory exists for the Armijo with damping rule, it may a good choice in applications, as it appears to possess the strengths of the other two methods.

\section*{Acknowledgments}
This work has been partly supported by the German Research Foundation (DFG) within the priority program SPP~1962 under contract number WE~6629/1-1 and Schu804/15-1, the Austrian Science Fund (FWF) within project no.~W1260-N35, and the German Academic Exchange Service (DAAD) within the Doctoral Programm 2017/18.
The authors are indebted to Laura Bittner for discussions about the model problem.

\bibliographystyle{plain}
\bibliography{references}
\end{document}

%% file: StochOpt-v2_shared.tex
\usepackage{lipsum}
\usepackage{amsfonts}
\usepackage{graphicx}
\usepackage{epstopdf}
\ifpdf
  \DeclareGraphicsExtensions{.eps,.pdf,.png,.jpg}
\else
  \DeclareGraphicsExtensions{.eps}
\fi

\usepackage{amsmath}
\usepackage{mathptmx}       
\usepackage{helvet}         
\usepackage{courier}        
\usepackage{overpic}
\usepackage{mathtools}
\usepackage{bbm}
\usepackage{caption,subcaption}
\usepackage{algpseudocode,algorithm,algorithmicx}

\usepackage{makeidx}

\usepackage{multicol}        
\usepackage[bottom]{footmisc}

\usepackage{subcaption}
\usepackage{cleveref}
\usepackage{stmaryrd}

\usepackage{tikz}
\usepackage{pgfplots}
\usetikzlibrary{spy,calc}

\newcommand{\pP}{\mathbb{P}}

\newcommand{\Go}{\partial D}
\newcommand{\Gi}{u}

\newcommand{\R}{\mathbb{R}}
\renewcommand{\d}{\textup{d}}
\newcommand{\dx}{\textup{d}x}
\newcommand{\ds}{\textup{d}s}

\newcommand{\id}{\operatorname{id}}

\newcommand{\LL}{{\mathcal{L}}}
\newcommand{\hH}{H_{\text{av}}^1(D)}
\newcommand{\hHh}{W_{\text{av},h}}
\newcommand{\EE}{\mathbb{E}}

\newtheorem{remark}{Remark}

\usepackage{amsopn}


%% file: zoom.tex
\newif\ifblackandwhitecycle
\gdef\patternnumber{0}

\pgfkeys{/tikz/.cd,
    zoombox paths/.style={
        draw=orange,
        very thick
    },
    black and white/.is choice,
    black and white/.default=static,
    black and white/static/.style={ 
        draw=white,   
        zoombox paths/.append style={
            draw=white,
            postaction={
                draw=black,
                loosely dashed
            }
        }
    },
    black and white/static/.code={
        \gdef\patternnumber{1}
    },
    black and white/cycle/.code={
        \blackandwhitecycletrue
        \gdef\patternnumber{1}
    },
    black and white pattern/.is choice,
    black and white pattern/0/.style={},
    black and white pattern/1/.style={    
            draw=white,
            postaction={
                draw=black,
                dash pattern=on 2pt off 2pt
            }
    },
    black and white pattern/2/.style={    
            draw=white,
            postaction={
                draw=black,
                dash pattern=on 4pt off 4pt
            }
    },
    black and white pattern/3/.style={    
            draw=white,
            postaction={
                draw=black,
                dash pattern=on 4pt off 4pt on 1pt off 4pt
            }
    },
    black and white pattern/4/.style={    
            draw=white,
            postaction={
                draw=black,
                dash pattern=on 4pt off 2pt on 2 pt off 2pt on 2 pt off 2pt
            }
    },
    zoomboxarray inner gap/.initial=5pt,
    zoomboxarray columns/.initial=2,
    zoomboxarray rows/.initial=2,
    subfigurename/.initial={},
    figurename/.initial={zoombox},
    zoomboxarray/.style={
        execute at begin picture={
            \begin{scope}[
                spy using outlines={%
                    zoombox paths,
                    width=\imagewidth / \pgfkeysvalueof{/tikz/zoomboxarray columns} - (\pgfkeysvalueof{/tikz/zoomboxarray columns} - 1) / \pgfkeysvalueof{/tikz/zoomboxarray columns} * \pgfkeysvalueof{/tikz/zoomboxarray inner gap} -\pgflinewidth,
                    height=\imageheight / \pgfkeysvalueof{/tikz/zoomboxarray rows} - (\pgfkeysvalueof{/tikz/zoomboxarray rows} - 1) / \pgfkeysvalueof{/tikz/zoomboxarray rows} * \pgfkeysvalueof{/tikz/zoomboxarray inner gap}-\pgflinewidth,
                    magnification=3,
                    every spy on node/.style={
                        zoombox paths
                    },
                    every spy in node/.style={
                        zoombox paths
                    }
                }
            ]
        },
        execute at end picture={
            \end{scope}
            \node at (image.north) [anchor=north,inner sep=0pt] {\subcaptionbox{\label{\pgfkeysvalueof{/tikz/figurename}-image}}{\phantomimage}};
            \node at (zoomboxes container.north) [anchor=north,inner sep=0pt] {\subcaptionbox{\label{\pgfkeysvalueof{/tikz/figurename}-zoom}}{\phantomimage}};
     \gdef\patternnumber{0}
        },
        spymargin/.initial=0.5em,
        zoomboxes xshift/.initial=1,
        zoomboxes right/.code=\pgfkeys{/tikz/zoomboxes xshift=1},
        zoomboxes left/.code=\pgfkeys{/tikz/zoomboxes xshift=-1},
        zoomboxes yshift/.initial=0,
        zoomboxes above/.code={
            \pgfkeys{/tikz/zoomboxes yshift=1},
            \pgfkeys{/tikz/zoomboxes xshift=0}
        },
        zoomboxes below/.code={
            \pgfkeys{/tikz/zoomboxes yshift=-1},
            \pgfkeys{/tikz/zoomboxes xshift=0}
        },
        caption margin/.initial=4ex,
    },
    adjust caption spacing/.code={},
    image container/.style={
        inner sep=0pt,
        at=(image.north),
        anchor=north,
        adjust caption spacing
    },
    zoomboxes container/.style={
        inner sep=0pt,
        at=(image.north),
        anchor=north,
        name=zoomboxes container,
        xshift=\pgfkeysvalueof{/tikz/zoomboxes xshift}*(\imagewidth+\pgfkeysvalueof{/tikz/spymargin}),
        yshift=\pgfkeysvalueof{/tikz/zoomboxes yshift}*(\imageheight+\pgfkeysvalueof{/tikz/spymargin}+\pgfkeysvalueof{/tikz/caption margin}),
        adjust caption spacing
    },
    calculate dimensions/.code={
        \pgfpointdiff{\pgfpointanchor{image}{south west} }{\pgfpointanchor{image}{north east} }
        \pgfgetlastxy{\imagewidth}{\imageheight}
        \global\let\imagewidth=\imagewidth
        \global\let\imageheight=\imageheight
        \gdef\columncount{1}
        \gdef\rowcount{1}
        
    },
    image node/.style={
        inner sep=0pt,
        name=image,
        anchor=south west,
        append after command={
            [calculate dimensions]
            node [image container,subfigurename=\pgfkeysvalueof{/tikz/figurename}-image] {\phantomimage}
            node [zoomboxes container,subfigurename=\pgfkeysvalueof{/tikz/figurename}-zoom] {\phantomimage}
        }
    },
    color code/.style={
        zoombox paths/.append style={draw=#1}
    },
    connect zoomboxes/.style={
    spy connection path={\draw[draw=none,zoombox paths] (tikzspyonnode) -- (tikzspyinnode);}
    },
    help grid code/.code={
        \begin{scope}[
                x={(image.south east)},
                y={(image.north west)},
                font=\footnotesize,
                help lines,
                overlay
            ]
            \foreach \x in {0,1,...,9} { 
                \draw(\x/10,0) -- (\x/10,1);
                \node [anchor=north] at (\x/10,0) {0.\x};
            }
            \foreach \y in {0,1,...,9} {
                \draw(0,\y/10) -- (1,\y/10);                        \node [anchor=east] at (0,\y/10) {0.\y};
            }
        \end{scope}    
    },
    help grid/.style={
        append after command={
            [help grid code]
        }
    },
}

\newcommand\phantomimage{%
    \phantom{%
        \rule{\imagewidth}{\imageheight}%
    }%
}
\newcommand\zoombox[2][]{
    \begin{scope}[zoombox paths]
        \pgfmathsetmacro\xpos{
            (\columncount-1)*(\imagewidth / \pgfkeysvalueof{/tikz/zoomboxarray columns} + \pgfkeysvalueof{/tikz/zoomboxarray inner gap} / \pgfkeysvalueof{/tikz/zoomboxarray columns} ) + \pgflinewidth
        }
        \pgfmathsetmacro\ypos{
            (\rowcount-1)*( \imageheight / \pgfkeysvalueof{/tikz/zoomboxarray rows} + \pgfkeysvalueof{/tikz/zoomboxarray inner gap} / \pgfkeysvalueof{/tikz/zoomboxarray rows} ) + 0.5*\pgflinewidth
        }
        \edef\dospy{\noexpand\spy [
            #1,
            zoombox paths/.append style={
                black and white pattern=\patternnumber
            },
            every spy on node/.append style={#1},
            x=\imagewidth,
            y=\imageheight
        ] on (#2) in node [anchor=north west] at ($(zoomboxes container.north west)+(\xpos pt,-\ypos pt)$);}
        \dospy
        \pgfmathtruncatemacro\pgfmathresult{ifthenelse(\columncount==\pgfkeysvalueof{/tikz/zoomboxarray columns},\rowcount+1,\rowcount)}
        \global\let\rowcount=\pgfmathresult
        \pgfmathtruncatemacro\pgfmathresult{ifthenelse(\columncount==\pgfkeysvalueof{/tikz/zoomboxarray columns},1,\columncount+1)}
        \global\let\columncount=\pgfmathresult
        \ifblackandwhitecycle
            \pgfmathtruncatemacro{\newpatternnumber}{\patternnumber+1}
            \global\edef\patternnumber{\newpatternnumber}
        \fi
    \end{scope}
}